# False discovery rate control with multivariate $p$-values

Zhiyi Chi *

*Department of Statistics*
*University of Connecticut*
*215 Glenbrook Road, U-4120*
*Storrs, CT 06269*
*e-mail:* zchi@stat.uconn.edu

**Abstract:** Multivariate statistics are often available as well as necessary in hypothesis tests. We study how to use such statistics to control not only false discovery rate (FDR) but also positive FDR (pFDR) with good power. We show that FDR can be controlled through nested regions of multivariate $p$-values of test statistics. If the distributions of the test statistics are known, then the regions can be constructed explicitly to achieve FDR control with maximum power among procedures satisfying certain conditions. On the other hand, our focus is where the distributions are only partially known. Under certain conditions, a type of nested regions are proposed and shown to attain (p)FDR control with asymptotically maximum power as the pFDR control level approaches its attainable limit. The procedure based on the nested regions is compared with those based on other nested regions that are easier to construct as well as those based on more straightforward combinations of the test statistics.

## 1. Introduction

In multiple hypothesis tests, it is common to evaluate nulls with univariate statistics. This especially has been the case for tests based on FDR control [2, 12, 13, 16, 20, 22, 23]. On the other hand, for hypotheses on high dimensional data, such as those in classification or recognition for complex signals, multivariate statistics in general are prerequisite for satisfactory results [1, 6, 24]. Such hypotheses each involves a sample of random vectors, from which a multivariate statistic is derived to capture critical features of the sample. Given the conceptual appeal of FDR control, it is natural to ask how it can be achieved using multivariate statistics.

The FDR of a multiple testing procedure is defined as $E[V/(R \vee 1)]$, where $R$ is the number of rejected nulls and $V$ that of rejected true nulls [2]. In addition to FDR, power and pFDR [20] are two important measures to assess the

* Research partially supported by NSF DMS 0706048 and NIH MH 68028. The author is thankful to the comments of the reviewers.

performance of a procedure. Recall that

$$\text{power} = E\left[\frac{R-V}{(n-N)\vee 1}\right], \quad \text{pFDR} = E[V/R \mid R > 0],$$

where $n$ is the number of nulls, and $N$ that of true nulls. The importance of power is well appreciated in the FDR literature [2, 13, 14, 20, 22]. In contrast, the issue of pFDR seems more subtle. Oftentimes, as follow-up actions can ensue only after *some* rejections are made, pFDR is more relevant than FDR. However, unlike FDR, in general pFDR is not necessarily controllable at a desirable level, say below 0.4. The reason is that in many cases, test statistics cannot provide strong enough evidence to assess the nulls, especially when the data distribution is only partially known and the number of observations for each null is small. The controllability of pFDR can strongly affect power. For the well-known Benjamini-Hochberg (BH) procedure [2], if its FDR control parameter is below the minimum attainable pFDR, then its power tends to 0 as $n \to \infty$ [7]. In light of this, power and pFDR should be considered together when designing testing procedures.

A direct way to improve power and pFDR control is to collect more observations for each null. However, this may not be feasible due to constraints on resources. On the other hand, if the observations can be viewed from different aspects each containing some unique information, then the aspects may be exploited together to yield more substantive evidence.

The approach of the paper is to first establish FDR control based on multivariate $p$-values, and then evaluate power and pFDR control. Among procedures that attain the same pFDR, the one with the highest power is preferred. Section 2 sets up notations and recalls known results. It then gives an example to illustrate when multivariate $p$-values may be useful for pFDR control. Section 3 presents a general FDR controlling procedure which uses an arbitrary family of nested regions in the domain of $p$-values. Then, it shows that if the data distribution under true nulls and that under false nulls are both known, then the nested regions can be chosen in such a way that the procedure has the maximum power among those with the same pFDR while satisfying certain consistency conditions. However, since full knowledge about data distributions is usually unavailable, the emphasis of the section is FDR control based on nested regions that approximate the optimal regions. Under certain conditions, the approximating regions are ellipsoids under an $L^\varepsilon$-norm, where $\varepsilon > 0$ in general is a non-integer.

Section 4 analyzes the power of the procedure based on the approximating regions. It shows that under certain conditions, the power is asymptotically maximized as the pFDR tends to the minimum attainable level. The procedure is compared with several others, including those that work "directly" on test statistics instead of $p$-values, for example, procedures that rejects nulls with large $L^a$-norms of the test statistics. It will be seen that only for $a < 0$, the "direct" procedures may attain the same pFDR as the procedure based on the approximating regions. The section also considers a procedure based on nested rectangle regions in the domain of $p$-values and shows that it has the same level of pFDR control as the one based on the approximating regions. Although less powerful, the procedure is simpler to compute.

Section 5 considers examples of $t$ and $F$ statistics. Section 6 reports a simulation study on the procedures considered in previous sections. Section 7 concludes with some remarks. Most of the technical details are collected in the Appendix.

## 2. Preliminaries

### 2.1. Notation

Denote by $K$ the dimension of a multivariate $p$-value. Points in $\mathbb{R}^K$ will be taken as column vectors. With a little abuse of notation, for $f\colon A \to \mathbb{R}$ with $A \subset \mathbb{R}^K$, $\sup f$ will denote the essential supremum of $f$, i.e. $\inf\{a\colon \ell(f^{-1}(a,\infty)) = 0\}$, where $\ell(\cdot)$ is the Lebesgue measure. If $\xi_{i1}, \ldots, \xi_{iK}$ are marginal or conditional $p$-values under a null $H_i$, then $\boldsymbol{\xi}_i = (\xi_{i1}, \ldots, \xi_{iK})'$ will be referred to as a multivariate $p$-value associated with $H_i$. The discussion is under a random effects model as follows [10, 13]. Denoting by $a \in (0,1)$ the proportion of false nulls and $\theta_i = \mathbf{1}\{H_i \text{ is false}\}$,

$$
\begin{aligned}
&(\theta_i, \boldsymbol{\xi}_i) \text{ are i.i.d. such that } \theta_i \sim \text{Bernoulli}(a) \text{ and} \\
&\text{given } \theta_i = 0, \ \xi_{i1}, \ldots, \xi_{iK} \text{ are i.i.d. } \sim \text{Unif}(0,1), \\
&\text{given } \theta_i = 1, \ \boldsymbol{\xi}_i \sim G \text{ with density } g.
\end{aligned} \quad (2.1)
$$

The density of $\boldsymbol{\xi}_i$ is then $1 - a + ag$. The assumption that $\xi_{i1}, \ldots, \xi_{iK}$ are independent under true $H_i$ should be checked carefully. Generally speaking, it should be problem-dependent to design test statistics with independent $p$-values [5]. One situation in which independence may arise is where multiple data sets on the same nulls are collected independently following different protocols, e.g., with different experiment designs being used or different physical attributes being recorded. In this situation, observations in different data sets may not follow the same distributions, and hence cannot be combined into larger i.i.d. samples. Nevertheless, the $p$-values derived respectively from them can be combined into multivariate $p$-values with independent coordinates.

Recall that, for univariate $p$-values $\xi_1, \ldots, \xi_n$, given FDR control parameter $\alpha \in (0,1)$, the BH procedure rejects $H_i$ with

$$\xi_i \leq \tau = \sup\left\{t \in [0,1] : \frac{t}{\alpha} \leq \frac{R(t) \vee 1}{n}\right\}, \quad R(t) = \#\left\{i : \xi_i \leq t\right\}.$$

Under the random effects model (2.1), the FDR actually realized by the BH procedure is $(1-a)\alpha$ [3, 11, 22], implying that the FDR can be arbitrarily small. On the other hand, the "local FDR" associated with each $H_i$ is [10]

$$P\left(\theta_i = 0 \mid \xi_i\right) = \frac{1-a}{1 - a + ag(\xi_i)} \geq \frac{1-a}{1 - a + a\sup g}.$$

It is not hard to see that the inequality applies to multivariate $p$-values as well. Then, following [8], the next result can be established.

**Proposition 2.1.** *Under* (2.1), *for any multiple testing procedure,*

$$\text{pFDR} \geq (1-a)\alpha_*, \quad \text{where} \ \ \alpha_* = \frac{1}{1 - a + a \sup g}. \tag{2.2}$$

Oftentimes, as $\sup g < \infty$, the pFDR is bounded from 0. In particular, if $a$ is small while $\sup g$ is only moderately large, the minimum attainable pFDR can be undesirably large. This is the basis of the next example.

### 2.2. An example

To further illustrate the role multivariate statistics may have for pFDR control, consider tests on $H_i : \boldsymbol{\mu}_i = (\mu_{X,i}, \mu_{Y,i}) = \mathbf{0}$ for $N(\boldsymbol{\mu}_i, \Sigma_i)$. Suppose for each $H_i$, a sample of $k$ i.i.d. $(X_{ij}, Y_{ij}) \sim N(\boldsymbol{\mu}_i, \Sigma_i)$ is collected. If it is known that $\Sigma_i \equiv \text{diag}(1,1)$ and under false $H_i$, $\mu_{X,i} = \mu_{Y,i} = 1$, then by Neyman-Pearson lemma, among procedures using fixed thresholding, the uniformly most powerful one is to reject $H_i$ if and only if $\bar{X}_i + \bar{Y}_i$ is greater than a suitable threshold value, where $\bar{X}_i = (1/k) \sum_j X_{ij}$ and likewise for $\bar{Y}_i$. That is to say, in this case univariate statistics are the best choice.

However, in most cases in practice, complete knowledge on data distributions is unavailable. If both $\Sigma_i$ and $\boldsymbol{\mu}_i$ under $H_i$ are unknown, then $\bar{X}_i + \bar{Y}_i$ cannot be used as test statistics and $t$ statistics are called for. Imagine a data analyst has computed the $t$ statistics of $\bar{X}_{ij}$ and those of $\bar{Y}_{ij}$ for each $H_i$, denoted $t_{X,i}$ and $t_{Y,i}$, respectively. While FDR control can be done with either $t_{X,i}$ or $t_{Y,i}$, the issue here is pFDR control.

Suppose the number of nulls is large and due to constraints on resources, $k = 9$ for each $H_i$. Suppose the data analyst knows that for each $H_i$, $\Sigma_i$ is diagonal and that for false $H_i$, $\mu_{X,i} > 0$ and $\mu_{Y,i} > 0$. However, he does not know that for false $H_i$, $\mu_{X,i}/\sigma_{X,i} = .5$ and $\mu_{Y,i}/\sigma_{Y,i} = .4$, where $\sigma_{X,i}^2$ and $\sigma_{Y,i}^2$ are the diagonal entries of $\Sigma_i$. If the fraction of false nulls is 5%, then, by using $t_{X,i}$ alone, the minimum attainable pFDR is $\approx .289$ and, by using $t_{Y,i}$ alone, the bound is even higher ($\approx .447$). The lower bounds are a consequence of Proposition 2.1. *No* procedure that only uses $t_{X,i}$ or $t_{Y,i}$ can get a pFDR lower than the bounds.

One way to attain lower pFDR is to increase $k$, which may require significantly more resources. When resources are limited, a sensible solution is to exploit both $t_{X,i}$ and $t_{Y,i}$ or, equivalently, their marginal $p$-values. This then raises the question of pFDR control using multivariate $p$-values.

## 3. FDR control using nested regions of *p*-values

### 3.1. General description

Let $\{D_t, \ 0 \leq t \leq 1\}$ be a family of Borel sets in $[0,1]^K$ such that

$$\begin{aligned} &D_1 = [0,1]^K, \ \ \ell(D_t) = t, \ \ D_s \subset D_t, \ \ 0 \leq s < t \leq 1 \\ &\{D_t\} \text{ is right-continuous, i.e., } D_t = \cap_{s>t} D_s, \ \ t \in [0,1). \end{aligned} \tag{3.1}$$

The most familar sets satifying (3.1) are perhaps $D_t = [0, t]$. Let $\boldsymbol{\xi}_1, \ldots, \boldsymbol{\xi}_n$ be the $p$-values associated with $H_1, \ldots, H_n$. Define

$$R(t) = \sum_{i=1}^n \mathbf{1}\{\boldsymbol{\xi}_i \in D_t\}, \quad V(t) = \sum_{i=1}^n (1-\theta_i)\mathbf{1}\{\boldsymbol{\xi}_i \in D_t\}.$$

**Description** Given FDR control parameter $\alpha \in (0,1)$,

$$\text{reject } H_i \text{ if and only if } \boldsymbol{\xi}_i \in D_\tau, \tag{3.2}$$

$$\text{where } \tau = \sup\left\{t \in [0,1] : \frac{t}{\alpha} \leq \frac{R(t) \vee 1}{n}\right\}. \qquad \square$$

**Theorem 3.1.** *For procedure* (3.2), $\text{FDR} = (1-a)\alpha$.

*Proof.* Since $D_t$ is right-continuous, $\boldsymbol{\xi}_i \in D_t \iff s_i \leq t$, where

$$s_i = \inf\{t \in [0,1] : \boldsymbol{\xi}_i \in D_t\}. \tag{3.3}$$

Therefore, procedure (3.2) rejects the same set of nulls as the BH procedure applied to $s_1, \ldots, s_n$ does. By $\ell(D_t) = t$, $s_i \sim \text{Unif}(0,1)$ under true $H_i$ and hence $s_i$ are univariate $p$-values. Theorem 3.1 then follows from [22]. $\square$

In general, a nested family of Borel sets in $[0,1]^K$ can often be parameterized so that procedure (3.2) is applicable to them.

**Proposition 3.1.** *Let* $\{\Gamma_u, u \in I\}$ *be a family of Borel sets in* $[0,1]^K$, *where* $I$ *is an interval in* $\mathbb{R}$, *such that* $\Gamma_u \subset \Gamma_v$ *for* $u < v$ *and* $\{\Gamma_u\}$ *is right-continuous. Suppose* $h(u) := \ell(\Gamma_u)$ *is continuous and strictly increasing with* $\inf h = 0$ *and* $\sup h = 1$. *For* $t \in (0,1)$, *define* $D_t = \Gamma_{h^{-1}(t)}$. *Also define* $D_0 = \cap \Gamma_u$ *and* $D_1 = [0,1]^K$. *Then procedure* (3.2) *based on* $D_t$ *attains* $\text{FDR} = (1-a)\alpha$.

As $\ell(D_t) \equiv t$, $D_t$ will be referred to as the regularization of $\Gamma_u$. Since nested regions naturally occur as decision regions in hypothesis tests, as seen below, by regularization, a test can turn into a FDR controlling procedure.

**Example 3.1.** (a) Suppose a test rejects a null if and only if $\min \xi_k \leq u$, where $\boldsymbol{\xi} = (\xi_1, \ldots, \xi_K)'$ is the associated $p$-value and $u$ a threshold value. The corresponding rejection region is $\Gamma_u = \{\boldsymbol{x} \in [0,1]^K : \min x_k \leq u\}$. Then $\{\Gamma_u, u \in [0,1]\}$ is an increasing family of sets. Since $h(u) = 1 - (1-u)^K$, procedure (3.2) applies to $D_t = \Gamma_{h^{-1}(t)}$ with $h^{-1}(t) = 1 - (1-t)^{1/K}$. Note that, in the Šidák procedure, when $K$ hypotheses are tested simultaneously, $h^{-1}(t)$ is the significance level for each hypothesis in order to attain familywise significance level $t$.

(b) Suppose a test rejects a null if and only if $\prod \xi_k \leq u$, where $u > 0$ is a threshold value. The corresponding rejection region is $\Gamma_u = \{\boldsymbol{x} \in [0,1]^K : \prod x_k \leq u\}$. For $K = 2$, $h(u) = u(1 + \ln u^{-1})$. In general, $h(u) = P(\prod U_k \leq u)$, with $U_k$ i.i.d. $\sim \text{Unif}(0,1)$. Since $-\ln U_k$ has density $e^{-x}\mathbf{1}\{x > 0\}$, $h(u) = 1 - F_K(-\ln u)$, where $F_K$ is the Gamma distribution with $K$ degrees of freedom and scale parameter 1. $\square$

For $\{D_t, 0 \le t \le 1\}$ having a regular representation, procedure (3.2) has an equivalent description more amenable to numerical evaluation. Suppose there is a function $0 \le J \le 1$, such that $D_t = \{\boldsymbol{x} \in [0,1]^K : J(\boldsymbol{x}) \le t\}$. It is easy to see $D_1 = [0,1]^K$ and $D_t$ is right-continuous for $t \in [0,1)$. Since $s_i = \inf\{t \in [0,1] : J(\boldsymbol{\xi}_i) \le t\} = J(\boldsymbol{\xi}_i)$, the next description obtains.

**Equivalent Description** Given FDR control level parameter $\alpha \in (0,1)$,

$$\text{apply the BH procedure to } s_i = J(\boldsymbol{\xi}_i). \tag{3.4}$$

That is, sort $s_1, \ldots, s_n$ into $s_{(1)} \le s_{(2)} \le \cdots \le s_{(n)}$. Define $s_{(0)} = 0$ and set $l = \max\{k \ge 0 : s_{(k)}/\alpha \le k/n\}$. Then reject $H_i$ if $s_i \le s_{(l)}$. □

EXAMPLE 3.1 (continued)

(a) Since $h(u) = 1 - (1-u)^K$ is strictly increasing, $h^{-1}$ and $D_t = \{\boldsymbol{x} \in [0,1]^K : \min x_k \le h^{-1}(t)\} = \{\boldsymbol{x} \in [0,1]^K : h(\min x_k) \le t\}$. Therefore $J(\boldsymbol{x}) = 1 - (1 - \min x_k)^K$.
(b) In this case $D_t = \{\boldsymbol{x} \in [0,1]^K : h(\prod x_k) \le t\}$, where $h(u) = 1 - F_K(-\ln u)$. Then $J(\boldsymbol{x}) = 1 - F_K(-\sum \ln x_k)$. For $K = 2$, since $h(u) = u(1 + \ln u^{-1})$, $J(x,y) = xy[1 - \ln x - \ln y]$.

### 3.2. Regions with maximum power under consistency condition

If the distribution under true nulls and that under false nulls are known, then, in light of Neyman-Pearson lemma, it is natural to ask if FDR can be controlled with maximum possible power using the likelihood ratios of the test statistics. Some works have been done on this idea [18, 21]. We next show that the idea is correct under certain conditions and can be realized by procedure (3.2) with an appropriate nested family $\{D_t, t \in [0,1]\} \subset [0,1]^K$.

Let $\boldsymbol{X} = (X_1, \ldots, X_K)' \in \mathbb{R}^K$ be a test statistic. Suppose that under true nulls, $\boldsymbol{X} \sim Q_0$ with density $q_0$ and under false nulls, $\boldsymbol{X} \sim Q_1$ with density $q_1$. Our construction of $D_t$ is based on a familiar transformation of $\boldsymbol{X}$ into multivariate p-values. Denote by $f_k(x_1, \ldots, x_k)$ the marginal density of $X_1, \ldots, X_k$ under $Q_0$. Clearly $q_0 = f_K$.

**Lemma 3.1.** Let $\phi(\boldsymbol{x}) = (\phi_1(\boldsymbol{x}), \ldots, \phi_K(\boldsymbol{x}))'$, $\boldsymbol{x} \in \mathbb{R}^K$, such that

$$\phi_1(\boldsymbol{x}) = Q_0(X_1 \le x_1), \quad \phi_k(\boldsymbol{x}) = Q_0(X_k \le x_k \mid X_s = x_s, \ s < k), \ k > 1.$$

Let $\boldsymbol{\xi} = \phi(\boldsymbol{X})$, i.e., $\xi_1$ is the p-value of $X_1$ and for $k > 1$, $\xi_k$ is the conditional p-value of $X_k$. Suppose i) $\text{sppt}(q_1) \subset \text{sppt}(q_0)$, where $\text{sppt}(q) := \{\boldsymbol{x} : q(\boldsymbol{x}) > 0\}$, ii) all $f_k$ are continuous and iii) $q_1$ is continuous on $\text{sppt}(q_0)$. Then 1) $\phi : \text{sppt}(q_0) \to [0,1]^K$ is continuous and 1-to-1; 2) under true nulls, $\xi_1, \ldots, \xi_K$ are i.i.d. $\sim \text{Unif}(0,1)$; and 3) under false nulls, $\boldsymbol{\xi}$ has a continuous density $g(\boldsymbol{x}) := q_1(\phi^{-1}(\boldsymbol{x}))/q_0(\phi^{-1}(\boldsymbol{x}))$ on $E := \phi(\text{sppt}(q_0))$, which is open with $\ell([0,1]\setminus E) = 1$.

Condition i) is not restrictive because hypothesis testing is trivial when $\boldsymbol{X} = \boldsymbol{x} \in \operatorname{sppt}(q_1) \setminus \operatorname{sppt}(q_0)$ (cf. [18]). Conditions ii) and iii) are used to make sure all the transformations involved on $\boldsymbol{X}$ are still well-defined random variables. Under these conditions, any rule based on the likelihood ratio of $\boldsymbol{X}$ has an equivalent based on $g(\boldsymbol{\xi})$. In the rest of the section, we will only consider tests on $p$-values.

A multiple testing procedure can be regarded as a deterministic or random function $\boldsymbol{\delta}$ that maps each data point to an $n$-tuple $(\delta_1, \ldots, \delta_n)$ with $\delta_i = \mathbf{1}\{H_i \text{ is rejected}\}$. In our setup, the data point is an $n$-tuple $(\boldsymbol{\xi}_1, \ldots, \boldsymbol{\xi}_n)$ jointly distributed with $\boldsymbol{\theta} = (\theta_1, \ldots, \theta_n)$. We need two conditions on $\boldsymbol{\delta}$.

(A) For $n \geq 1$, $\boldsymbol{\delta}$ and $\boldsymbol{\theta}$ are independent conditional on $\boldsymbol{\xi}_1, \ldots, \boldsymbol{\xi}_n$, i.e.,

$$P(\boldsymbol{\delta} = \boldsymbol{a} \,|\, \boldsymbol{\xi}_1, \ldots, \boldsymbol{\xi}_n, \boldsymbol{\theta} = \boldsymbol{b}) = P(\boldsymbol{\delta} = \boldsymbol{a} \,|\, \boldsymbol{\xi}_1, \ldots, \boldsymbol{\xi}_n), \quad \boldsymbol{a}, \boldsymbol{b} \in \{0, 1\}^n.$$

Most multiple testing procedures are deterministic functions of test statistics and therefore satisfy condition (A). The condition means that the observed test statistics contain all the available information on $\boldsymbol{\theta}$; if there is any prior knowledge on $\boldsymbol{\theta}$, it has already been fully incorporated into $\boldsymbol{\xi}$ and hence any randomness introduced into $\boldsymbol{\delta}$ is a "pure guess".

The second condition imposes some consistency on $\boldsymbol{\delta}$. For an $n$-tuple $S = (\boldsymbol{x}_1, \ldots, \boldsymbol{x}_n)$ with $\boldsymbol{x}_i \in [0,1]^K$, denote $R(S; \boldsymbol{\delta}) = \sum_{i=1}^n \delta_i(S)$ and $\hat{F}(\boldsymbol{x}; S)$ the empirical distribution function

$$\hat{F}(\boldsymbol{x}; S) = \#\{i : x_{ik} \leq x_k, \ k = 1, \ldots, K\}/n.$$

(B) For any sequence of $n_k$-tuples $S_k$ with $n_k \to \infty$, if $\hat{F}(\boldsymbol{x}; S_k)$ converges in the sense that $\sup_{\boldsymbol{x}} |\hat{F}(\boldsymbol{x}; S_k) - F(\boldsymbol{x})| \to 0$ for a distribution function $F$, then $R(S_k; \boldsymbol{\delta})/n_k$ converges in probability.

Basically, the condition requires that, when $\boldsymbol{\delta}$ is applied to samples with similar empirical distributions, it should reject similar fractions of nulls from them. Loosely speaking, that means as far as the fraction of rejected nulls is concerned, $\boldsymbol{\delta}$ has to "stick to" a single way of testing, rather than alternate between different ways for different data sets.

For $0 \leq u < \infty$, define

$$\Gamma_u = \{\boldsymbol{x} \in [0, 1]^K : g(\boldsymbol{x}) \geq u\}. \tag{3.5}$$

Although $\{\Gamma_u\}$ is decreasing instead of increasing, its regularization can be made increasing. Let $h(u) = \ell(\Gamma_u)$. Then $h$ is decreasing, $h(0) = 1$ and $h(u) \to 0$ as $u \to \infty$. Define $h^*(t) = \inf\{u \geq 0 : h(u) \leq t\}$.

**Proposition 3.2.** *Suppose conditions i)–iii) in Lemma 3.1 are satisfied and $h$ is continuous. Then $D_t = \Gamma_{h^*(t)}$ satisfies (3.1). Let $\alpha \in (0,1)$. Then, among procedures that satisfy conditions (A) and (B) and attain $\mathrm{FDR} \leq (1-a)\alpha$, procedure (3.2) with $D_t$ belongs to those which asymptotically have the maximum power as $n \to \infty$. Furthermore, the following statements on the procedure hold:*

1) it always rejects the same set of nulls as the BH procedure does when applied to $p_i = h(g(\boldsymbol{\xi}_i))$, $i = 1, \ldots, n$; 2) for $\alpha > \alpha_*$, the power is asymptotically positive and $\mathrm{pFDR} \to \mathrm{FDR} = (1-a)\alpha$; and 3) for $\alpha < \alpha_*$, the power is asymptotically 0 and $\mathrm{pFDR} \to (1-a)\alpha_*$.

**Example 3.2.** To illustrate that in general condition (B) is needed in Proposition 3.2, let $K = 1$. First consider the case where the $p$-values $\xi_1, \ldots, \xi_n$ are i.i.d. $\sim F(t) = (1-a)t + aG(t) \in C^{(1)}([0,1])$ such that $t/F(t)$ is strictly increasing and $G(t)$ is linear on $[t_1, t_2]$, where $0 < t_1 < t_2 < 1$. It is easy to see that $t/F(t)$ is strictly concave on $[t_1, t_2]$. Given $c \in (0, 1)$, consider the following randomized procedure. Draw $U \sim \mathrm{Unif}(0,1)$. If $U > c$, reject and only reject nulls with $\xi_i \in [0, t_1]$; otherwise, reject and only reject nulls with $\xi_i \in [0, t_2]$. As $n \to \infty$, the empirical distribution of $\xi_i$ converges to $F$. However, conditional on $U > c$, $R_n/n \to F(t_1)$, while conditional on $U < c$, $R_n/n \to F(t_2)$. Therefore, the procedure satisfies condition (A) but not (B). It can be seen that $\mathrm{pFDR} \to (1-a)\alpha$ with $\alpha = (1-c)t_1/F(t_1) + ct_2/F(t_2)$ and power $\to (1-c)G(t_1) + cG(t_2) = G(t_c)$, with $t_c = (1-c)t_1 + ct_2$.

Consider the BH procedure when it is applied to $\xi_1, \ldots, \xi_n$ with control parameter $\alpha$. Since $1/F'(0) < \alpha < 1$, by [12], the procedure asymptotically has power $G(t^*)$, where $t^* \in (0, 1)$ such that $t^*/F(t^*) = \alpha$. Since $t/F(t)$ is strictly increasing, $t_1 < t^* < t_2$. On the other hand, since $t/F(t)$ is strictly concave on $[t_1, t_2]$, $t_c/F(t_c) > \alpha$. As a result, $t_c > t^*$. Therefore, asymptotically, although the BH procedure has the same pFDR level as the randomized procedure, it is strictly less powerful.

Finally, given $c \in (0, 1)$, by small variation to $G$ on $[t_1, t_2]$, one can construct $G$ which is smooth and strictly concave, such that the above conclusions still hold. By Proposition 3.2, the most powerful procedure (3.2) satisfying condition (B) in this case is the BH procedure and hence is strictly less powerful than the randomized procedure at the same pFDR level. □

As noted in a discussion in [4], by either accepting all nulls with probability $1 - \alpha$ or rejecting all of them with probability $\alpha$, it is guaranteed that $\mathrm{FDR} \leq \alpha$; however, the FDR attained in this way is useless, because it cannot say how well one can learn from *the* data being analyzed. Without some coherence of a procedure, one can hardly make a sensible evaluation of its performance in a particular instance based on a measure defined as a long term average, as the measure incorporates not only the way of testing chosen for the data at hand, but also others that are potentially very different. The same comments apply to pFDR as well. Condition (B) aims to impose some coherence, which is possible if the data follows the law of large numbers. This is similar to the ergodicity assumption, whereby long term average can be approximated by an average over a single large sample.

The construction in this section requires full knowledge of the density $g$, which is often unavailable. However, if $g$ is known to possess some regularities, then it is possible to apply procedure (3.2) to regular shaped $D_t$ with reasonable power. This possibility is explored next.

### *3.3. Nonconstant approximation of lowest order*

In many cases, the following is true for the distribution $G$ of $p$-values under false nulls:

$$G \text{ has density } g \in C([0,1]^K) \text{ with } g(\boldsymbol{x}) < g(\boldsymbol{0}) < \infty \text{ for all } \boldsymbol{x} \neq \boldsymbol{0}. \qquad (3.6)$$

Under (3.6), smaller $p$-values are stronger evidence against nulls, nevertheless the strength is bounded. By Proposition 2.1, the minimum attainable pFDR is $(1-a)\alpha_*$, where $\alpha_*$ is now equal to $1/[1-a+ag(\boldsymbol{0})]$. Following Taylor's expansion, suppose for some $\gamma_k > 0$ and $\varepsilon > 0$,

$$g(\boldsymbol{x}) = g(\boldsymbol{0})(1 - \boldsymbol{\gamma}'\boldsymbol{x}^\varepsilon + r(\boldsymbol{x})) \text{ with } r(\boldsymbol{x}) = o(|\boldsymbol{x}|^\varepsilon) \text{ as } |\boldsymbol{x}| \to 0, \qquad (3.7)$$

where $\boldsymbol{x}^\varepsilon$ denotes $(x_1^\varepsilon, \ldots, x_K^\varepsilon)'$. It is perhaps desirable and expected to be true that $\varepsilon$ is a positive integer. However, under regular conditions, this usually is not the case. As will be seen in Section 5, for the upper-tail $p$-values associated with $t$ or $F$ statistics, $\varepsilon$ usually is a fraction of 1. More generally, $g(\boldsymbol{x}) = g(\boldsymbol{0})(1 - \sum \gamma_k x_k^{\varepsilon_k}) + o(\sum x_k^{\varepsilon_k})$ with $\varepsilon_k > 0$ possibly different. However, for simplicity, this case will not be discussed.

Rewrite the region in (3.5) as $\{\boldsymbol{x} \in [0,1]^K : g(\boldsymbol{x}) \geq g(\boldsymbol{0})(1-u)\}$. For $0 < u \ll 1$, the region is approximately $\{\boldsymbol{x} \in [0,1]^K : \boldsymbol{\gamma}'\boldsymbol{x}^\varepsilon \leq u\}$, suggesting that the latter may be used in procedure (3.4) with reasonable power. In general, for $\boldsymbol{\nu} = (\nu_1, \ldots, \nu_K)'$ with $\nu_k \geq 0$ and $\sum \nu_k > 0$, define $\Gamma_u$ as

$$\Gamma_u(\boldsymbol{\nu}) = \{\boldsymbol{x} \in [0,1]^K : \boldsymbol{\nu}'\boldsymbol{x}^\varepsilon \leq u\}. \qquad (3.8)$$

Then by procedure (3.4), the following procedure obtains.

**Control based on regions (3.8)** Given FDR control parameter $\alpha \in (0,1)$,

$$\text{apply the BH procedure to } s_i = h(\boldsymbol{\nu}'\boldsymbol{\xi}_i^\varepsilon; \boldsymbol{\nu}), \text{ where} \qquad (3.9)$$

$$h(u; \boldsymbol{\nu}) = \ell(\Gamma_u(\boldsymbol{\nu})) = \int_0^1 \cdots \int_0^1 \mathbf{1}\left\{\sum \nu_k x_k^\varepsilon \leq u\right\} d\boldsymbol{x}. \qquad \square$$

Procedure (3.9) is "scale invariant" in $\boldsymbol{\nu}$, i.e., the set of nulls rejected by using $\Gamma_u(c\boldsymbol{\nu})$ is the same for $c > 0$. If $K = 1$, then the procedure is simply the BH procedure and the parameter $\boldsymbol{\nu} = \nu_1$ has no effect on its performance. However, when $K > 1$, the power of procedure (3.9) depends on $\boldsymbol{\nu}$. To analyze the power, denote

$$V_\varepsilon = \int_0^1 \cdots \int_0^1 \mathbf{1}\left\{\sum x_k^\varepsilon \leq 1\right\} d\boldsymbol{x}.$$

The next lemma will be used.

**Lemma 3.2.** *Given $\boldsymbol{\nu} = (\nu_1, \ldots, \nu_K)'$ with $\nu_k > 0$, let $\bar{\nu}$ be the geometric mean of $\nu_k$, i.e., $\bar{\nu} = (\nu_1 \cdots \nu_K)^{1/K}$. Then for $0 < u \le \min \nu_k$,*

$$h(u; \boldsymbol{\nu}) = V_\varepsilon \left(\frac{u}{\bar{\nu}}\right)^{K/\varepsilon}, \tag{3.10}$$

$$\int_{\Gamma_u(\boldsymbol{\nu})} x_k^\varepsilon \, d\boldsymbol{x} = \frac{V_\varepsilon}{K+\varepsilon} \frac{\bar{\nu}}{\nu_k} \left(\frac{u}{\bar{\nu}}\right)^{K/\varepsilon+1}, \quad k = 1, \ldots, K. \tag{3.11}$$

*Furthermore, as $u \downarrow 0$,*

$$\int_{\Gamma_u(\boldsymbol{\nu})} g = g(\mathbf{0}) V_\varepsilon \left(1 - \frac{u}{K+\varepsilon} \sum \frac{\gamma_k}{\nu_k}\right) \left(\frac{u}{\bar{\nu}}\right)^{K/\varepsilon} + o(u^{K/\varepsilon+1}). \tag{3.12}$$

### 3.4. Special cases

In most cases, $h(u; \boldsymbol{\nu})$ is complicated to evaluate. There are two cases that allow tractable numerical evaluation of $h(u; \boldsymbol{\nu})$. The first case is $K = 2$. Suppose $\nu_2 \ge \nu_1 > 0$. For $u \ge \nu_1 + \nu_2$, it is clear $h(u; \boldsymbol{\nu}) = 1$. For $0 < u < \nu_1 + \nu_2$,

$$\begin{aligned} h(u; \boldsymbol{\nu}) = & \left[\left(\frac{u - \nu^*}{\nu_*}\right) \vee 0\right]^{1/\varepsilon} + \frac{\Gamma(1/\varepsilon)^2}{2\varepsilon \Gamma(2/\varepsilon)} \left(\frac{u}{\bar{\nu}}\right)^{2/\varepsilon} \\ & \times \left[F\left(\frac{\nu_*}{u}; \frac{1}{\varepsilon}, 1 + \frac{1}{\varepsilon}\right) - F\left(1 - \frac{\nu^*}{u}; \frac{1}{\varepsilon}, 1 + \frac{1}{\varepsilon}\right)\right] \end{aligned} \tag{3.13}$$

where $\nu_* = \min(\nu_1, \nu_2)$, $\nu^* = \max(\nu_1, \nu_2)$ and $F(x; a, b)$ is the Beta distribution with parameters $a$ and $b$. See Appendix A.1 for a proof of (3.13).

The second case is $\nu_1 = \cdots = \nu_K$ and $\varepsilon = 1$, where $h(u; \boldsymbol{\nu})$ can be evaluated by recursion. Due to scale invariance, let $\nu_k = 1$. Then

$$h_K(u) := h(u; 1, \ldots, 1) = \int_0^1 \cdots \int_0^1 \mathbf{1}\left\{\sum x_k \le u\right\} d\boldsymbol{x}$$

is piecewise polynomial, such that $h_K(u) = h_{K, \lfloor u \rfloor}(\{u\})$, where $\lfloor u \rfloor$ is the largest integer no greater than $u$, $\{u\} = u - \lfloor u \rfloor$, and

$$h_{K,i}(t) := h_K(t+i) = \sum_{k=0}^{K} A_K(i,k) t^k, \quad t \in [0, 1). \tag{3.14}$$

Since $h_{K,K}(t) \equiv 1$, $h_{K,0}(t) = t^K/K!$ and, for $i = 1, \ldots, K-1$,

$$\begin{aligned} h_{K,i}(t) &= \int_0^t h_{K-1}(t-x) \, dx + \int_t^1 h_{K-1}(t-x) \, dx \\ &= \int_0^t h_{K-1,i}(t-x) \, dx + \int_t^1 h_{K-1,i-1}(1+t-x) \, dx \\ &= \int_0^t h_{K-1,i}(x) \, dx + \int_t^1 h_{K-1,i-1}(x). \end{aligned}$$

It follows that for $i = 1, \ldots, K-1$ and $k = 1, \ldots, K$,

$$\begin{cases} A_K(0,0) = \cdots = A_K(0, K-1) = 0, \quad A_K(0,K) = 1/K! \\ A_K(i,0) = \sum_{k=0}^{K-1} A_{K-1}(i-1,k)/(k+1), \\ A_K(i,k) = [A_{K-1}(i,k-1) - A_{K-1}(i-1,k-1)]/k, \\ A_K(K,0) = 1, \quad A_K(K,k) = 0. \end{cases} \quad (3.15)$$

By $h_1(u) = u$, the initial conditions are $A_1(0,0) = 0$, $A_1(0,1) = 1$. These relations together with (3.14) can be used to compute $h_K(u)$.

## 4. Analysis of power

Recall that power $= E[\frac{R-V}{(n-N)\vee 1}]$, with $n$ the number of nulls and $N$ that of true nulls. As our focus is the case $n \gg 1$, we shall consider the limit $\text{Pow}(\alpha)$ of power at pFDR level $(1-a)\alpha$ as $n \to \infty$. In general, closed form formulas for $\text{Pow}(\alpha)$ are not available. To get a handle on $\text{Pow}(\alpha)$, our approach is to look at how fast it drops to 0 as $\alpha \downarrow \alpha_*$, by approximating $\text{Pow}(\alpha)$ as a linear combination of $(\alpha - \alpha_*)^a$, $a \in [0, \infty)$, or for that matter, $G(D_t)$ as a linear combination of $t^a$, with $D_t$ a family of nested regions used by procedure (3.2). Thus, the analysis is essentially a type of Taylor's expansion, which can provide useful qualitative information for comparing powers of different procedures.

Our analysis will only yield approximations of low orders. It remains to be seen how high order approximations can be obtained. In order to apply the results in section 3, we shall assume

$$g \text{ satisfies (3.6) and (3.7) such that } \ell(\Gamma_u) \text{ is continuous}$$

where $\Gamma_u$ is defined in (3.5).

### 4.1. Dependence on parameter values

For $\alpha$ close to $\alpha_*$, the dependency of the power of procedure (3.9) on $\boldsymbol{\nu}$ can be characterized as follows.

**Proposition 4.1.** *Fix $\nu_k > 0$. Then for procedure (3.9), the minimum attainable pFDR is $(1-a)\alpha_*$ and*

$$\text{Pow}(\alpha) \sim g(\mathbf{0}) V_\varepsilon \left( \frac{K+\varepsilon}{a\alpha_*^2 g(\mathbf{0})} \bigg/ \sum \frac{\bar{\nu}\gamma_k}{\nu_k} \right)^{K/\varepsilon} (\alpha - \alpha_*)^{K/\varepsilon} \quad \text{as } \alpha \downarrow \alpha_*. \quad (4.1)$$

Due to the scale invariance of procedure (3.9), let $\nu_1 \cdots \nu_K = \gamma_1 \cdots \gamma_K$. Let $\lambda_k = \gamma_k/\nu_k$. Then $\sum \bar{\nu}\gamma_k/\nu_k = \bar{\gamma} \sum \lambda_k$, which is minimized under the constraint $\lambda_1 \cdots \lambda_K = 1$ if and only if $\lambda_k \equiv 1$. It follows that as $\alpha \downarrow \alpha_*$, $\text{Pow}(\alpha)$ asymptotically is maximized if $\boldsymbol{\nu} = \boldsymbol{\gamma}$ and

$$\sup_{\boldsymbol{\nu}} \text{Pow}(\alpha) \sim g(\mathbf{0}) V_\varepsilon \left[ \frac{1 + \varepsilon/K}{a\alpha_*^2 g(\mathbf{0})\bar{\gamma}} \right]^{K/\varepsilon} (\alpha - \alpha_*)^{K/\varepsilon}.$$

Recall that procedure (3.2) based on the regions in (3.5) has the maximum power for $\alpha > \alpha_*$ among procedures satisfying conditions (A) and (B). The next result says that, at $\boldsymbol{\nu} = \boldsymbol{\gamma}$, procedure (3.9) and this procedure are asymptotically equivalent, i.e., as $\alpha \downarrow \alpha_*$, they not only have about the same power, but also reject about the same sets of nulls.

**Proposition 4.2.** *For $\alpha > \alpha_*$, let $\mathrm{Pow}_o(\alpha)$ be the limit of power of procedure (3.2) based on (3.5), and $\mathrm{Pow}(\alpha)$ that of procedure (3.9) with $\boldsymbol{\nu} = \boldsymbol{\gamma}$. Then $\mathrm{Pow}_o(\alpha)/\mathrm{Pow}(\alpha) \to 1$, as $\alpha \downarrow \alpha_*$.*

*Moreover, let $\mathcal{V}_o$ and $\mathcal{D}_o$ be the sets of true nulls and false nulls, respectively, that are rejected by the first procedure and $\mathcal{V}$ and $\mathcal{D}$ those by the second one. Let $r_D(\alpha)$ be the in-probability limit of $|\mathcal{D}_o \triangle \mathcal{D}|/|\mathcal{D}_o \cap \mathcal{D}|$ and $r_V(\alpha)$ that of $|\mathcal{V}_o \triangle \mathcal{V}|/|\mathcal{V}_o \cap \mathcal{V}|$ as $n \to \infty$. Then $r_D(\alpha) \to 0$ and $r_V(\alpha) \to 0$ as $\alpha \downarrow \alpha_*$.*

### 4.2. Other types of nested regions

In order to compare the power of procedure (3.9) and that of procedure (3.2) based on other types of nested regions, the following comparison lemma will be used, which says that if a nested family of regions can "round up" more false nulls, then procedure (3.2) based on the regions has more power.

**Lemma 4.1.** *Let $\{D_{it}\}$, $i = 1, 2$, be two families of Borel sets satisfying (3.1). Suppose $G(D_{it})$ are continuous in $t$ and there is $T \in (0, 1)$, such that for $0 < t \leq T$, $G(D_{1t}) < G(D_{2t})$. Given $\alpha \in (\alpha_*, 1)$, for the procedure based on $D_{it}$, let $\tau_i$ be defined as in (3.2). Assume that as $n \to \infty$, $\tau_i \xrightarrow{P} t_i^* \in (0, T)$. Then $\mathrm{Pow}_1(\alpha) < \mathrm{Pow}_2(\alpha)$.*

To start with, for $\Gamma_u(\boldsymbol{\nu})$ in (3.8), by (3.10), the regularization is $D_t = \Gamma_{u(t)}(\boldsymbol{\nu})$ with $u(t) = \bar{\nu}(t/V_\varepsilon)^{\varepsilon/K}$ if $0 < t \ll 1$. Then by (3.12),

$$G(D_t) = g(\mathbf{0}) \left[1 - \frac{(t/V_\varepsilon)^{\varepsilon/K} \bar{\nu}}{K + \varepsilon} \sum \frac{\gamma_k}{\nu_k}\right] t + o(t^{\varepsilon/K+1}), \quad \text{as } t \to 0. \qquad (4.2)$$

**Example 4.1.** A common rule is to reject a null if and only if $z = \prod \xi_k$ is small. The rejection regions are $\Gamma'_u = \{\boldsymbol{x} \in [0,1]^K : \prod x_k \leq u\}$. We next show that for $\alpha \approx \alpha_*$, procedure (3.2) based on $\Gamma'_u$ has strictly less power than procedure (3.9) for any $\boldsymbol{\nu}$ with $\nu_1 \cdots \nu_K > 0$. Roughly, the reason is that, for $u \ll 1$, most of $\Gamma'_u$ is spread around the boundary surfaces $x_k = 0$ where the density of false nulls is lower than that around $\mathbf{0}$.

First consider $K = 2$. By Example 3.1(b), the regularization of $\Gamma'_u$ is $D_{1t} = \{(x, y) \in [0, 1]^2 : xy \leq h^{-1}(t)\}$ with $h(u) = u(1 + \ln u^{-1})$. Denote the regularization for $\Gamma_u(\boldsymbol{\nu})$ by $D_{2t}$. As $\alpha \downarrow \alpha_*$, $t_i^* \to 0$. Thus, by Lemma 4.1, in order to compare the powers for $\alpha \approx \alpha_*$, it is enough to compare $G(D_{1t})$ and $G(D_{2t})$ for $t \ll 1$. Recall $\int_{D_{it}} dx\, dy = t$. Fix $a_i \in (0, \gamma_i)$ and $\eta > 0$ such that

$g(x, y) < g(\mathbf{0})(1 - a_1 x^\varepsilon - a_2 y^\varepsilon)$ for $x, y \in [0, \eta]$. By $G(D_{1t}) = \int_{D_{1t}} g$,

$$G(D_{1t}) \leq g(\mathbf{0}) \left[ \int_{D_{1t}} dx\, dy - \int_{\substack{xy \leq h^{-1}(t) \\ 0 < x, y < \eta}} (a_1 x^\varepsilon + a_2 y^\varepsilon)\, dx\, dy \right]$$
$$= g(\mathbf{0})[t - (a_1 + a_2)\eta^\varepsilon h^{-1}(t)/\varepsilon] + O(h^{-1}(t)^{1+\varepsilon}) \quad \text{as } t \to 0.$$

By (4.2), $G(D_{2t}) = g(\mathbf{0})[t - Ct^{1+\varepsilon/2}] + o(t^{1+\varepsilon/2})$ with $C > 0$ a constant. Since $h^{-1}(t)/t^{\varepsilon/2+1} \to \infty$ as $t \to 0$, $G(D_{1t}) < G(D_{2t})$ for $t \ll 1$. Thus by Lemma 4.1, $\text{Pow}_1(\alpha) < \text{Pow}_2(\alpha)$ for $\alpha \approx \alpha_*$. Furthermore, the next result implies $\text{Pow}_1(\alpha) = o(\text{Pow}_2(\alpha))$ as $\alpha \downarrow \alpha_*$.

**Proposition 4.3.** *Under the setup in Lemma 4.1, let $D_{2t}$ be the regularization of the regions $\Gamma_u(\boldsymbol{\nu})$ in (3.8). If $t_1^* \to 0$ as $\alpha \to \alpha_*$ and*

$$D_{1t} \downarrow \{\mathbf{0}\} \quad \text{and} \quad \frac{g(\mathbf{0})t - G(D_{1t})}{g(\mathbf{0})t - G(D_{2t})} \to M \in [1, \infty] \quad \text{as } t \to 0,$$

*then $\text{Pow}_1(\alpha)/\text{Pow}_2(\alpha) \to (1/M)^{K/\varepsilon}$ as $\alpha \downarrow \alpha_*$.*

The case $K > 2$ can be treated likewise; see Appendix A.2. □

**Example 4.2.** Normal quantile transformations of $p$-values have been used as a convenient representation of data in multiple testing [9]. Denote $\bar\Phi(x) = P(Z \geq x)$ with $Z \sim N(0,1)$. Let $w_k > 0$, $k = 1, \ldots, K$. Consider the rule that rejects a null if any only if $Q(\boldsymbol{\xi}) = \sum_k w_k \bar\Phi^{-1}(\xi_k)$ is large. The corresponding rejection regions are $\Gamma'_u = \{\boldsymbol{x} \in [0,1]^K : Q(\boldsymbol{x}) \geq u\}$. As $\bar\Phi(\xi) \sim N(0,1)$ for $\xi \sim \text{Unif}(0,1)$,

$$\ell(\Gamma'_u) = P\left(\sum_{k=1}^K w_k \bar\Phi^{-1}(\xi_k) \geq u\right) \qquad \xi_i \text{ i.i.d. } \sim \text{Unif}(0,1)$$
$$= P\left(\sum_{k=1}^K w_k Z_k \geq u\right) = \bar\Phi(u/w) \qquad Z_i \text{ i.i.d. } \sim N(0,1),$$

where $w = \sqrt{\sum_k w_k^2}$. Then the regularization of $\Gamma'_u$ is $D_{1t} = \{\boldsymbol{x} \in [0,1]^K : Q(\boldsymbol{x}) \geq w\bar\Phi^{-1}(t)\}$. In Appendix A.2, it is shown that

$$\overline{\lim_{t \downarrow 0}} \frac{\ln[g(\mathbf{0})t - G(D_{1t})]}{\ln t} \leq b := \left[1 - \frac{\varepsilon}{(1+\varepsilon)K}\right]^{-1}. \tag{4.3}$$

On the other hand, with $D_{2t}$ as in Example 4.1,

$$\lim_{t \downarrow 0} \frac{\ln[g(\mathbf{0})t - G(D_{2t})]}{\ln t} = 1 + \varepsilon/K.$$

For $K > 1$, $b < 1 + \varepsilon/K$. As a result, $\frac{g(\mathbf{0})t - G(D_{1t})}{g(\mathbf{0})t - G(D_{2t})} \to \infty$ as $t \downarrow 0$. It then follows from Proposition 4.3 that $\text{Pow}_1(\alpha) = o(\text{Pow}_2(\alpha))$ as $\alpha \downarrow \alpha_*$. □

**Example 4.3** (Nested rectangle regions). An easy way to get nested regions is as follows. Let $0 \leq f_k(t) \leq 1$ be nondecreasing continuous functions on $[0, 1]$, such that $\prod f_k(t) = t$ and $f_k(t) \to 0$ as $t \to 0$. Let

$$D_{1t} = [0, f_1(t)] \times \cdots \times [0, f_K(t)] = \{\boldsymbol{x} \in [0,1]^K : f_k^*(x_k) \leq t, \text{ all } k\}, \quad (4.4)$$

where $f_k^*(x) = \sup\{u : f_k(u) < x\}$. Now procedure (3.4) is the BH procedure applied to $s_i = J(\boldsymbol{\xi}_i) := \max_k f_k^*(\xi_{ik})$. By Appendix A.2,

$$G(D_{1t}) = g(\boldsymbol{0})\left(1 - \frac{1}{1+\varepsilon}\sum \gamma_k f_k(t)^\varepsilon\right)t + o\left(t \sum f_k(t)^\varepsilon\right) \quad \text{as } t \downarrow 0. \quad (4.5)$$

Among all $f_k$ with $\prod f_k(t) \equiv t$, $(\bar{\gamma}/\gamma_k)^{1/\varepsilon}t^{1/K}$ minimize $\sum \gamma_k f_k(t)^\varepsilon$. Thus the maximum of $G(D_{1t})$ is $g(\boldsymbol{0})[1 - \bar{\gamma}Kt^{\varepsilon/K}/(1+\varepsilon)]t + o(t^{1+\varepsilon/K})$.

Let $D_{2t}$ be the regularization for $\Gamma_u(\boldsymbol{\nu})$. For $K = 1$, $D_{1t} = D_{2t}$. For $K > 1$, by (4.2), $G(D_{2t})$ is asymptotically maximized with value $g(\boldsymbol{0})[1 - \bar{\gamma}KLt^{\varepsilon/K}]t$ if $\boldsymbol{\nu} = \boldsymbol{\gamma}$, where $L = (1/V_\varepsilon)^{\varepsilon/K}/(K + \varepsilon)$. By Appendix A.3, $L < 1/(1+\varepsilon)$ for $K \geq 2$. Thus, for $\alpha \approx \alpha_*$, the maximum power of procedure (3.2) based on $D_{1t}$ is strictly less than the one based on $D_{2t}$. By Proposition 4.3,

$$\frac{\text{Pow}_1(\alpha)}{\text{Pow}_2(\alpha)} \to \frac{1}{V_\varepsilon}\left(\frac{1+\varepsilon}{K+\varepsilon}\right)^{K/\varepsilon} \in (0, 1), \quad \text{as } \alpha \downarrow \alpha_*.$$

Therefore, the power by using the rectangles is of the same order as that by using $\Gamma_u(\boldsymbol{\nu})$, albeit lower.

Procedure (3.2) based on rectangles has some advantages, even though the power is not maximum. First, rectangles are much easier to construct than the regularization of $\Gamma_u(\boldsymbol{\nu})$. Second, for $g$ satisfying (3.7), there is no need to know $\varepsilon$. Indeed, it suffices to try $f_k(t) = c_k t^{1/K}$ with $\prod c_k = 1$. By (3.4), the next procedure obtains.

**Control based on rectangles** Given FDR control parameter $\alpha \in (0, 1)$,

$$\text{apply the BH procedure to } s_i = \left[\max_k \frac{\xi_{ik}}{c_k}\right]^K, \quad (4.6)$$

where $c_k > 0$ satisfy $c_1 \cdots c_K = 1$. □

By Appendix A.2, for procedure (4.6),

$$\text{Pow}(\alpha) \sim g(\boldsymbol{0})\left(\frac{1+\varepsilon}{a\alpha_*^2 g(\boldsymbol{0})} \Big/ \sum \gamma_k c_k^\varepsilon\right)^{K/\varepsilon}(\alpha - \alpha_*)^{K/\varepsilon} \quad \text{as } \alpha \downarrow \alpha_*. \quad (4.7)$$

In particular, if $c_k = (\bar{\gamma}/\gamma_k)^{1/\varepsilon}$, the power is asymptotically maximized.

Now suppose $\xi_{ik}$ are upper-tail probabilities of test statistics $X_{ik}$ and $f_k(t) = t^{1/K}$. Then procedure (4.6) rejects $H_i$ if and only if $p_{ik} \leq \tau^{1/K}$ for all $k$, where $\tau$ is random. If $\tau$ is small, then procedure (4.6) may be viewed as one that rejects $H_i$ with large $\min_k X_{ik}$, which makes it seem unnatural as $\max_k X_{ik}$ is more often used in testing. However, it will be seen in Example 4.4 that procedures based on $\max_k X_{ik}$ in general cannot attain the same level of pFDR control as procedure (4.6).

### 4.3. Direct combination of test statistics

Procedure (3.9) requires $p$-values of test statistics. Oftentimes, procedures that only use simple combinations of test statistics seem more desirable because they do not have to evaluate $p$-values. However, as seen next, in many cases such procedures cannot attain pFDR control levels as low as procedure (3.9), and hence have strictly less power at low pFDR control levels.

We only consider test statistics $\boldsymbol{X}_i = (X_{i1}, \ldots, X_{iK})$ with $X_{ik}$ being independent under false $H_i$ as well as under true $H_i$. Let $X_{ik}$ follow $F_{0k}$ under true $H_i$ and $F_k$ under false $H_i$ and suppose $F_{0k}$ and $F_k$ have continuous densities $f_{0k}$ and $f_k$ respectively, such that $(0, \infty) \subset \mathrm{sppt}(f_k) \subset \mathrm{sppt}(f_{0k})$. Denote $r_k(x) = f_k(x)/f_{0k}(x)$ with $0/0$ set to 0, and suppose

$$r_k(x) \text{ is strictly increasing on } \mathrm{sppt}(f_k) \text{ but } \rho_k := \lim_{x \to \infty} r_k(x) < \infty, \quad (4.8)$$

i.e., for each $X_{ik}$, larger values are stronger evidence against $H_i$, however, the strength is bounded. Notice that, under (4.8), $\rho_k > \int r_k(x) f_{0k}(x)\,dx = \int f_k(x)\,dx = 1$.

Let $\xi_{ik}$ be the upper-tail $p$-value of $X_{ik}$. Then $\xi_{ik} = \overline{F}_{0k}(X_{ik}) = 1 - F_{0k}(X_{ik})$. Under false $H_i$, $\xi_{ik} \sim G_k(u) = 1 - F_k(F_{0k}^{-1}(1-u))$ with density

$$g_k(u) = \frac{f_k(\phi_k(u))}{f_{0k}(\phi_k(u))}, \quad \text{where} \quad \phi_k(u) = F_{0k}^{-1}(1-u). \quad (4.9)$$

Then $\boldsymbol{\xi}_i \sim G(\boldsymbol{x}) = \prod G_k(x_k)$ with density $g(\boldsymbol{x}) = \prod g_k(x_k)$. By (4.8) and (4.9), $g_k$ are strictly decreasing and $g_k(0) = \rho_k \in (1, \infty)$. Then $g(\boldsymbol{0}) = \prod \rho_k = \sup g$ and (3.6) is satisfied. If

$$g_k(u) = g_k(0)[1 - \gamma_k u^\varepsilon + o(u^\varepsilon)], \text{ as } u \downarrow 0, \ k = 1, \ldots, K, \quad (4.10)$$

then $g(\boldsymbol{x})$ satisfies (3.7) and so by Propositions 4.1, the minimum attainable pFDR level for procedure (3.9) is $(1-a)\alpha_*$, with $\alpha_* = 1/(1-a+ag(\boldsymbol{0}))$. In the following examples, we shall assume (4.10) holds.

**Example 4.4.** One common combination of $X_{ik}$ is $M_i = \max_k X_{ik}$. Under true $H_i$, $M_i \sim \prod F_{0k}$ and, under false $H_i$, $M_i \sim \prod F_k$. By only using $M_i$, the minimum attainable pFDR is $(1-a)/(1-a+aL)$, where

$$L = \sup_x \frac{\sum_k f_k(x) \prod_{j \neq k} F_j(x)}{\sum_k f_{0k}(x) \prod_{j \neq k} F_{0j}(x)}.$$

In Appendix A.2, it is shown that $L < g(\boldsymbol{0})$. Thus $(1-a)/(1-a+aL) > (1-a)\alpha_*$.

For procedures that reject $H_i$ if and only if $M_i$ is large enough, say, $M_i \geq T$, where $T \gg 1$ is a fixed threshold value, the minimum attainable pFDR can be even higher. Indeed, now $\mathrm{pFDR} \geq (1-a)/(1-a+aL')$, where

$$L' = \sup_{x \geq T} \frac{1 - \prod_k F_k(x)}{1 - \prod_k F_{0k}(x)} = (1 + \eta_T) \sup_{x \geq T} \frac{\sum_k \overline{F}_k(x)}{\sum_k \overline{F}_{0k}(x)}$$
$$\leq (1 + \eta_T) \max_k \rho_k, \qquad \text{with } \eta_T \to 0 \text{ as } T \to \infty.$$

Since $\max_k \rho_k$ can be much smaller than $g(\mathbf{0}) = \prod_k \rho_k$, the minimum pFDR can be significantly higher than $(1-a)\alpha_*$. □

**Example 4.5.** Another common combination of $X_{ik}$ is $M_i = \sum_k c_k X_{ik}$, where $c_k > 0$, such that $H_i$ is rejected if and only if $M_i$ is large enough. Without loss of generality, let $K = 2$ and $c_k = 1$. Suppose $X_{ik} \geq 0$, such that $f_{0k}(x) \sim x^{-s_k}$ and $f_k(x) \sim \rho_k x^{-s_k}$ as $x \to \infty$ with $s_k > 1$. Then $g(\mathbf{0}) = \rho_1 \rho_2$. Under true $H_i$, the density of $M_i$ is $\int_0^x f_{01}(t) f_{02}(x-t)\, dt$ and, under false $H_i$, it is $\int_0^x f_1(t) f_2(x-t)\, dt$. Then the minimum attainable pFDR by only using $M_i$ is $(1-a)/(1-a+aL)$, where

$$L = \sup_{x \geq T} r(x), \quad r(x) := \frac{\int_0^x f_1(t) f_2(x-t)\, dt}{\int_0^x f_{01}(t) f_{02}(x-t)\, dt},$$

with $T > 0$ a threshold value. In Appendix A.2, it is shown that $L < g(\mathbf{0})$. As a result, the minimum attainable pFDR is strictly greater than $(1-a)\alpha_*$.

The same conclusion holds if $M_i$ is a weighted $L^q$-norm of $\mathbf{X}_i$ with $q > 0$, i.e., $M_i = (\sum_k c_k X_{ik}^q)^{1/q}$. To see this, let $c_k = 1$. Under true $H_i$, $X_{ik}^q \sim G_{0k}(x) = F_{0k}(x^{1/q})$, and under false $H_i$, $X_{ik}^q \sim G_k(x) = F_k(x^{1/q})$. Then $G'_{0k}(x) \sim (1/q) x^{-t_k}$ and $G'_k(x) \sim (\rho_k/q) x^{-t_k}$, where $t_k = (s_k - 1)/q + 1 > 1$. Then the argument for $\sum_k X_{ik}$ can be applied to $G_{0k}$ and $G_k$. □

**Example 4.6.** Under the same conditions as in Example 4.5, assume further that $f_k(x)/f_{0k}(x) = \rho_k[1 - D_k x^{-d_k} + o(x^{-d_k})]$, with $D_k \neq 0$ and $d_k > 0$. Note that $D_k > 0$. It follows that as $x_1, \ldots, x_K \to \infty$,

$$\prod \frac{f_k(x_k)}{f_{0k}(x_k)} = \left[1 - (1+o(1)) \sum D_k x_k^{-d_k}\right] \times \prod \rho_k.$$

Therefore, in order to attain pFDR $= (1-a)\alpha$ with $\alpha \approx \alpha_*$, a null $H_i$ should be rejected if and only if $v_i := \sum D_k X_{ik}^{-d_k}$ is small. In particular, if $d_k = d$, then $v_i^{-1/d}$ is a weighted $L^{-d}$-norm of $\mathbf{X}_i$ and hence $H_i$ is rejected if and only if the norm is large.

Indeed, letting $\xi_{ik}$ be the upper $p$-value of $X_{ik}$, by the derivations in next section, $v_i = \sum[1 + o_k(\xi_{ik})] \nu_k \xi_{ik}^{\varepsilon_k}$, where $\nu_k = D_k(s_k - 1)^{\varepsilon_k}$, $\varepsilon_k = d_k/(s_k - 1)$ and $o_k(u) \to 0$ as $u \to 0$. Consequently, if $\varepsilon_k \equiv \varepsilon$, then the above procedure can be formulated as one based on $\boldsymbol{\xi}_i$, such that the associated rejection regions are approximately $\Gamma_u(\boldsymbol{\nu})$ in (3.9) when contracting to $\mathbf{0}$. Provided the regularizations of the rejection regions are readily available, they can be used in procedure (3.2) for FDR control as well, with approximately the same power as procedure (3.9) as $\alpha \downarrow \alpha_*$. □

## 5. Examples of special distributions

We next show $t$ and $F$ distributions satisfy (3.6) and (3.7). To this end, some general formulas are needed. Suppose $X_1, \ldots, X_K$ are test statistics that are independent not only under true nulls but also under false nulls. Let $F_{0k}$ and $f_{0k}$ be the marginal distribution and density of $X_k$ under a true null, and $F_k$ and $f_k$ those under a false null.

Suppose that, for some positive constants $a, b, c_k, d_k$ and $r_k$, $1 \le k \le K$,

$$f_{0k}(x) \sim a c_k x^{-a-1} \tag{5.1}$$

$$\frac{f_k(x)}{f_{0k}(x)} = r_k - d_k x^{-b} + o(x^{-b}), \quad \text{as } x \to \infty. \tag{5.2}$$

Then, as $u \downarrow 0$, $\psi_k(u) := F_{0k}^{-1}(1-u) \sim (c_k/u)^{1/a}$. Let $\xi_k$ be the upper-tail $p$-value of $X_k$ and $g_k(u)$ its density. By (4.9) and (5.2), as $u \downarrow 0$,

$$g_k(u) = r_k - d_k \psi_k(u)^{-b} + o(\psi_k(u)^{-b}) = r_k - c_k^{-b/a} d_k u^{b/a} + o(u^{b/a}), \tag{5.3}$$

implying $g_k(0) = r_k$. By independence, $\boldsymbol{\xi}$ has density $g(\boldsymbol{u}) = \prod g_k(u_k)$ which satisfies (3.6) and (3.7). Moreover, the parameters in (3.7) are

$$\varepsilon = b/a, \quad g(\boldsymbol{0}) = r_1 \cdots r_K, \quad \gamma_k = c_k^{-\varepsilon} d_k / r_k. \tag{5.4}$$

Note that for $F_{0k} = N(0, \sigma^2)$ and $F_k = N(\mu, \sigma^2)$, where $\mu > 0$ and $\sigma$ is known, (5.3) does not hold, as $\sup g_k = \sup f_k / f_{0k} = \infty$. For simplicity, we next only consider $K = 1$ and omit the index $k$.

### 5.1. t distribution

Let $F = t_{p,\delta}$, the noncentral $t$ distribution with $p$ dfs and noncentrality parameter $\delta \ge 0$ and $F_0 = t_p = t_{p,0}$. The density of $t_{p,\delta}$ is

$$t_{p,\delta}(x) = \frac{A e^{-\delta^2/2}}{(p+x^2)^{(p+1)/2}} \sum_{k=0}^{\infty} C_k (\sqrt{2}\delta)^k \left(\frac{x^2}{p+x^2}\right)^{k/2},$$

where $A = \frac{p^{p/2} \Gamma((p+1)/2)}{\sqrt{\pi} \Gamma(p/2)}, \quad C_k = \frac{1}{k!} \Gamma\left(\frac{p+1+k}{2}\right) \Big/ \Gamma\left(\frac{p+1}{2}\right).$

By $C_0 = 1$, $t_p(x) = A C_0 (p+x^2)^{-(p+1)/2} \sim A x^{-p-1}$ as $x \to \infty$. Let $z = 1 - x/\sqrt{p+x^2}$. Then $z \sim (p/2) x^{-2}$ and

$$\frac{t_{p,\delta}(x)}{t_p(x)} = e^{-\delta^2/2} \sum_{k=0}^{\infty} C_k (\sqrt{2}\delta)^k (1-z)^k$$

$$= e^{-\delta^2/2} \left[\sum_{k=0}^{\infty} C_k (\sqrt{2}\delta)^k - \frac{p}{2x^2} \sum_{k=1}^{\infty} k C_k (\sqrt{2}\delta)^k\right] + o(z).$$

Consequently, in (5.1) and (5.2), $a = p$, $b = 2$, $c = A/p$,

$$d = \frac{e^{-\delta^2/2} p}{2} \sum_{k=1}^{\infty} k C_k (\sqrt{2}\delta)^k < \infty, \quad r = e^{-\delta^2/2} \sum_{k=0}^{\infty} C_k (\sqrt{2}\delta)^k < \infty.$$

Then by (5.4),

$$\varepsilon = \frac{b}{a} = \frac{2}{p}, \quad \gamma = \frac{1}{2}\left[\frac{p\sqrt{\pi}\,\Gamma(p/2)}{\Gamma((p+1)/2)}\right]^{2/p} \sum_{k=1}^{\infty} k C_k(\sqrt{2}\delta)^k \Big/ \sum_{k=0}^{\infty} C_k(\sqrt{2}\delta)^k.$$

Clearly, for $p > 2$, $\varepsilon$ is a fraction of 1.

### 5.2. F distribution

Let $F = F_{p,q,\delta}$, the noncentral $F$ distribution with $(p, q)$ dfs and noncentrality parameter $\delta \geq 0$ and $F_0 = F_{p,q} = F_{p,q,0}$. Letting $\rho = p/q$ and $z = 1/(1 + \rho x)$, $F_{p,q,\delta}$ has density

$$f_{p,q,\delta}(x) = \frac{e^{-\delta/2}}{x}(1-z)^{p/2}z^{q/2}\sum_{k=0}^{\infty}\frac{(\delta/2)^k(1-z)^k}{k!B(p/2+k,q/2)},$$

where $B(x, y) = \Gamma(x)\Gamma(y)/\Gamma(x+y)$. The density of $F_{p,q}$ is $f_{p,q}(x) = f_{p,q,0}(x)$. As $x \to \infty$, $z \sim (1/\rho)x^{-1}$. It follows that

$$f_{p,q}(x) = \frac{1}{B(p/2,q/2)x}(1-z)^{p/2}z^{q/2} \sim \frac{\rho^{-q/2}}{B(p/2,q/2)}x^{-q/2-1},$$

$$\frac{f_{p,q,\delta}(x)}{f_{p,q}(x)} = e^{-\delta/2}\sum_{k=0}^{\infty}\frac{C_{p,q,k}}{k!}\left(\frac{\delta}{2}\right)^k(1-z)^k$$

$$= e^{-\delta/2}\left[\sum_{k=0}^{\infty}\frac{C_{p,q,k}}{k!}\left(\frac{\delta}{2}\right)^k - \frac{1}{\rho x}\sum_{k=1}^{\infty}\frac{C_{p,q,k}}{(k-1)!}\left(\frac{\delta}{2}\right)^k\right] + o(z),$$

where $C_{p,q,k} = \dfrac{B(p/2,q/2)}{B(p/2+k,q/2)}.$

Thus in (5.1) and (5.2),

$$a = q/2, \quad b = 1, \quad c = \frac{2\rho^{-q/2}}{qB(p/2,q/2)} = \frac{\rho^{-q/2}}{B(p/2,q/2+1)}$$

$$d = \frac{e^{-\delta/2}}{\rho}\sum_{k=1}^{\infty}\frac{C_{p,q,k}}{(k-1)!}\left(\frac{\delta}{2}\right)^k < \infty, \quad r = e^{-\delta/2}\sum_{k=0}^{\infty}\frac{C_{p,q,k}}{k!}\left(\frac{\delta}{2}\right)^k < \infty$$

and by (5.4), $\varepsilon = b/a = 2/q$,

$$\gamma = \left[B\left(\frac{p}{2}, \frac{q}{2}+1\right)\right]^{2/q}\sum_{k=1}^{\infty}\frac{C_{p,q,k}(\delta/2)^k}{(k-1)!}\Big/\sum_{k=0}^{\infty}\frac{C_{p,q,k}(\delta/2)^k}{k!}.$$

Clearly, for $q > 2$, $\varepsilon$ is a fraction of 1.

## 6. Numerical study

We report a simulation study on the procedures described in previous sections. The simulations are implemented in R language [17]. We focus on testing mean vectors of multivariate normals with only partial knowledge about their variances. Model (2.1) is used to sample $H_1, \ldots, H_n$, such that $H_i$ is "$\boldsymbol{\mu}_i = \mathbf{0}$ in $N(\boldsymbol{\mu}_i, \Sigma_i)$" and $\mathbf{1}\{H_i \text{ is true}\} \sim \text{Bernoulli}(a)$ with $a = .05$ or $.02$, where $\boldsymbol{\mu}_i \in \mathbb{R}^K$. We mimic the situation where the only available information is 1) under true $H_i$, $\Sigma_i$ is diagonal and 2) under false $H_i$, all the coordinates of $\boldsymbol{\mu}_i$ are positive. Since $\Sigma_i$ are not necessarily the same for different $H_i$ and there is no knowledge on their relations whatsoever, $\Sigma_i$ cannot be estimated by pooling the observations. Since $\Sigma_i$ are unknown and the values of $\boldsymbol{\mu}_i$ under false $H_i$ are also unknown, the tests have to rely on $t$-statistics. In the simulations, for each $H_i$, an i.i.d. sample $\boldsymbol{X}_{i1}, \ldots, \boldsymbol{X}_{i,\text{df}+1} \sim N(\boldsymbol{\mu}_i, \Sigma_i)$ is drawn and the test statistic is computed as $(T_{i1}, \ldots, T_{iK})'$, where $T_{ik} = \sqrt{\text{df}+1}\bar{X}_{ik}/S_{ik}$, with $\bar{X}_{ik}$ and $S_{ik}$ the mean and standard deviation of the $k$-th coordinates of $\boldsymbol{X}_{ij}$. It follows that for true $H_i$, $T_{ik} \sim t_n$ and for false $H_i$, $T_{ik} \sim t_{\text{df}, \sqrt{\text{df}+1}\mu_{ik}/\sigma_{ik}}$, where $\mu_{ik}$ is the $k$-th coordinate of $\boldsymbol{\mu}_i$ and $\sigma_{ik}$ the $k$-th diagonal entry of $\Sigma_i$. The corresponding $p$-value is $\boldsymbol{\xi}_i = (\xi_{i1}, \ldots, \xi_{iK})'$, with $\xi_{ik}$ the marginal upper-tail $p$-values of $T_{ik}$ under $t_n$.

### 6.1. Bivariate p-values

In this part, we compare procedures (3.9) and (4.6). Throughout, $K = 2$ and FDR control parameter $\alpha = .15$. For true $H_i$, $N(\boldsymbol{\mu}_i, \Sigma_i) = N(\mathbf{0}, I_K)$ and for false $H_i$, $N(\boldsymbol{\mu}_i, \Sigma_i) = N(\boldsymbol{\mu}, \Sigma(r))$, where the diagonal entries of $\Sigma(r)$ are equal to 1 and off-diagonal entries equal to $2r/(1+r^2)$. For $r = 0$, the coordinates of $\boldsymbol{\xi}_i$ are independent under false $H_i$. To examine the effect of dependency between the coordinates of $\boldsymbol{\xi}_i$ under false $H_i$, we also simulate with $r = \pm 1/5$. Each simulation makes 2000 runs, each run tests 5000 nulls. The (p)FDR and power are computed as Monte Carlo averages of the runs.

We conduct 3 groups of simulations, corresponding to $(\boldsymbol{\mu}, n) = (.6, .2, 8)$, $(.5, .5, 8)$ and $(2, 2, 2)$, respectively (Table 1). In each group, procedures (3.9) and (4.6) are implemented for 6 pairs of $a$ and $\Sigma(r)$, with $a = .05, .02$ and $r = 0, \pm 1/5$. For the pairs with $r = 0$, $\alpha_* = 1/[1 - a + ag(\mathbf{0})]$ is calculated using (5.4) and the results in Section 5.1. The value of $\alpha_*$ as well as those of $\varepsilon = 2/n$, $\gamma_1$ and $\gamma_2$ are given in Table 1. As is seen, for all the pairs $\alpha = .15 > \alpha_*$, so pFDR $= (1-a)\alpha$ is attainable by the procedures. On the other hand, if only the $i$th ($i = 1, 2$) coordinate of the $p$-value is used, pFDR $= (1-a)\alpha$ is not attainable because in this case the minimum attainable pFDR is $(1-a)\alpha_{*i}$ with $\alpha_{*i} > .15$.

TABLE 1
Parameters of the simulations in Section 6.1. The values of $\alpha_*$, $\gamma$ and $\varepsilon = 2/\text{df}$ are computed for the case where the coordinates of $p$-value are independent

|   | $\boldsymbol{\mu}$, df | $\alpha_*, \alpha_{*1}, \alpha_{*2}$ ($a = .05, .02$) | | $\boldsymbol{\gamma}$ | $\varepsilon$ |
|---|---|---|---|---|---|
| 1 | (.6, .2), 8 | $9.37 \times 10^{-3}, .18, .47$ | $2.31 \times 10^{-2}, .36, .69$ | (5.82, 3.51) | 1/4 |
| 2 | (.5, .5), 8 | $9.05 \times 10^{-3}, .30, .30$ | $2.23 \times 10^{-2}, .52, .52$ | (46.81, 46.81) | 1/4 |
| 3 | (2, 2), 2 | $2.88 \times 10^{-2}, .44, .44$ | $6.90 \times 10^{-2}, .67, .67$ | (27.69, 27.69) | 1 |

For $r = 0$ and $\alpha \approx \alpha_*$, by Proposition 4.1 and (4.7), procedures (3.9) and (4.6) approximately reach their respective maximum power at pFDR level $(1-a)\alpha$ if $\nu_k = \gamma_k$ in (3.9) and $c_k = c_k^{(0)} := (\bar\gamma/\gamma_k)^{1/\varepsilon}$ in (4.6). To see how the powers depend on $\nu_k$ and $c_k$ for $\alpha = .15$, the procedures are tested with

$$\begin{cases} \boldsymbol{\nu} = (\nu_1, \nu_2) = (s^\varepsilon \gamma_1, \gamma_2/s^\varepsilon) & \text{for procedure (3.9)}, \\ \boldsymbol{c} = (c_1, c_2) = (sc_1^{(0)}, c_2^{(0)}/s) & \text{for procedure (4.6)}, \end{cases} \quad (6.1)$$

where $s > 0$ is a tuning parameter. The reason why $s^\varepsilon$ instead of $s$ is used for $\boldsymbol{\nu}$ will be seen later. For $r = \pm 1/5$, the procedures are tested with the same sets of values $\boldsymbol{\nu}$ and $\boldsymbol{c}$ as well. For groups 1 and 2, (3.13) is used to calculate $h(u; \boldsymbol{\nu})$ for procedure (3.9). For group 3, as $\varepsilon = 1$, (3.15) is used.

The plots of power and (p)FDR *vs* $\log_2 s$ are shown in Figures 1–3 and labeled with "e" and "r" for procedures (3.9) and (4.6), respectively. The label "e" refers to "ellipsoid", due to the similarity of the nested regions in procedure (3.9) to Euclidean ellipsoids. For most of the plots, $a = .05$. The results for $a = .02$ are qualitatively the same, except that the power is lower and the pFDR is harder to control. For illustration, Figure 1 includes the plots of power and (p)FDR for $(\boldsymbol{\mu}, \mathrm{df}) = (.6, .2, 8)$, $r = 0$ and $a = .02$.

The results show that for $\alpha$ significantly greater than $\alpha_*$, the power may still exhibit patterns similar to that for $\alpha \approx \alpha_*$. First, by Proposition 4.1 and (4.7), for $\alpha \approx \alpha_*$, the maximum power of procedure (3.9) is strictly greater than that of (4.7). The left panels of Figures 1–3 show that this remains to be the case for $\alpha = .15$. Second, for $\boldsymbol{\nu}$ and $\boldsymbol{c}$ as in (6.1), as $\alpha \approx \alpha_*$, for both procedures, the power is approximately proportional to $(s^\varepsilon + 1/s^\varepsilon)^{-K/\varepsilon}$. As a result, the power curves of the procedures should be approximately symmetric, decreasing in $|\log_2 s|$ and parallel to each other. This holds quite well for $\alpha = .15$, except for the plots for procedure (4.6) in Figure 1, which exhibit moderate asymmetry that may be attributable to the unequal marginal distributions of the coordinates of the $p$-values.

In (6.1), it is necessary to use $s^\varepsilon$ to tune $\boldsymbol{\nu}$ in order to get a power curve parallel to the one for procedure (4.6). For $t$ distributions with df $> 2$, $\varepsilon < 1$, suggesting that procedure (3.9) is more sensitive to the change in $\boldsymbol{\nu}$ than procedure (4.6) is to the change in $\boldsymbol{c}$. However, since $\varepsilon$ is known, as the results show, the sensitivity is easy to address. For df $= 2$, $\varepsilon = 1$ and hence the power of procedure (3.9) is uniformly greater than that of procedure (4.6) (Figure 3).

The results also demonstrate the difference between pFDR and FDR. In the simulations, the FDR remains constant. However, as power decreases, the pFDR increases, sometimes quite rapidly. Unless power is high enough, the pFDR is strictly greater than the FDR. Note that this observation is made when the number of tested nulls is 5000. In theory, if power is positive, then pFDR $\to$ FDR as the number of nulls tends to $\infty$. The observed discrepancy between the pFDR and FDR is due to the fact that the number of nulls is not large enough for the asymptotic to take effect.

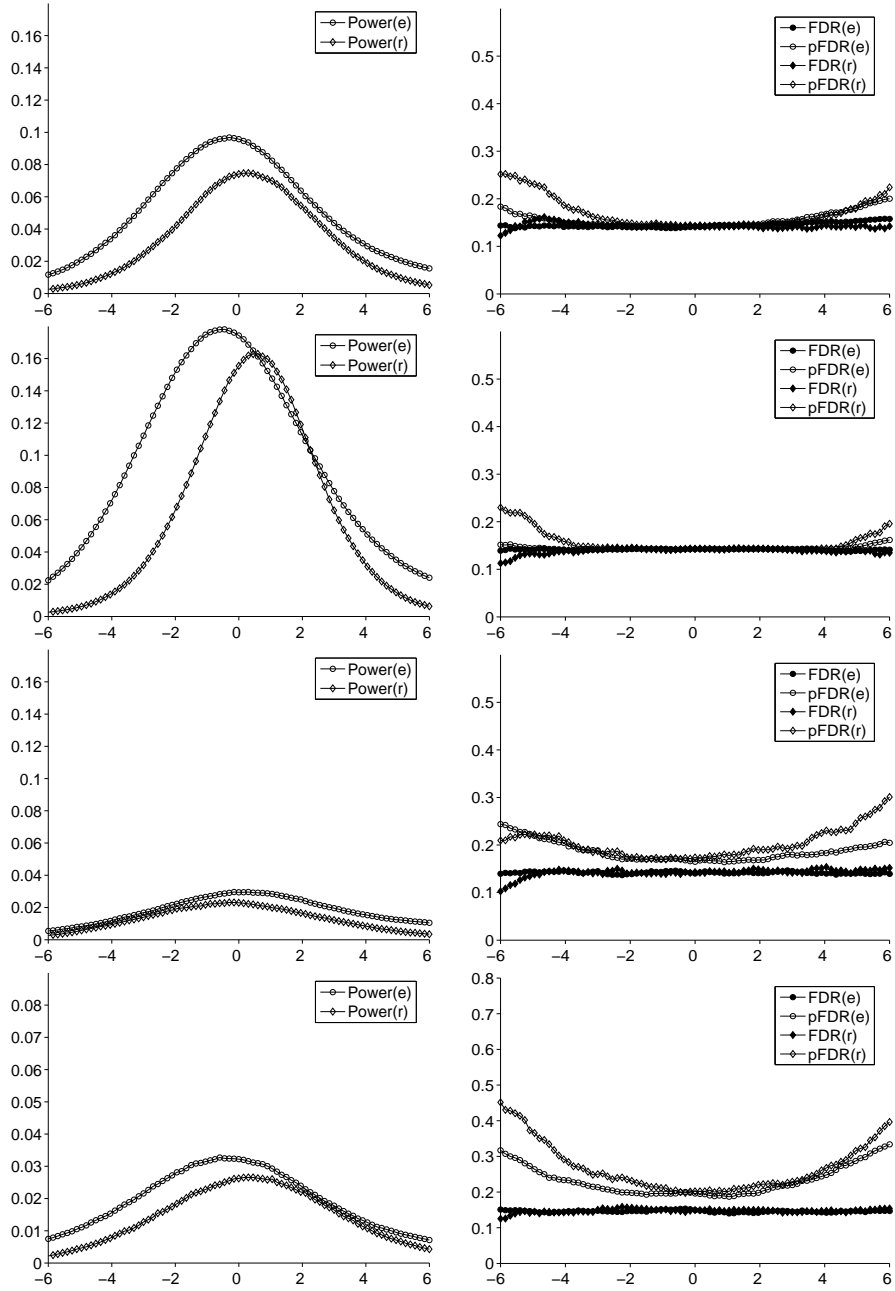

FIG 1. Power and (p)FDR vs $\log s$ for procedures (3.9) and (4.6): group 1 (cf. Section 6.1). Rows 1-3, $a = .05$, $r = 0, 1/2, -1/2$. Row 4, $a = .02$, $r = 0$.

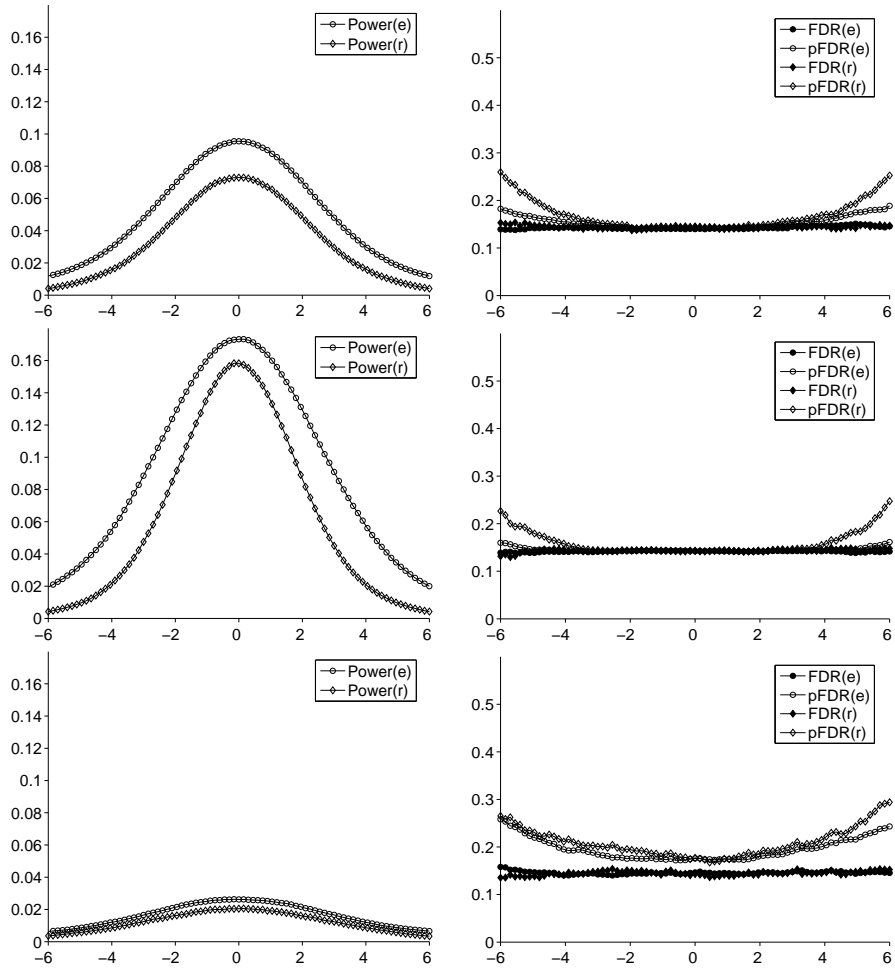

Fig 2. Power and (p)FDR *vs* $\log s$ for procedures (3.9) and (4.6): group 2 (cf. Section 6.1). $a = .05$, $r = 0, 1/2, -1/2$.

Finally, as seen from Figures 1 and 2, statistical dependency between the coordinates of the test statistics may have significant influence on power and pFDR control. Nevertheless, the modality and symmetry of the power curves are quite stable. Furthermore, the effects of correlations are not obvious in Figure 3, where $\varepsilon = 1$ and the joint distribution of the *p*-values is symmetric. This apparent stability remains to be explained.

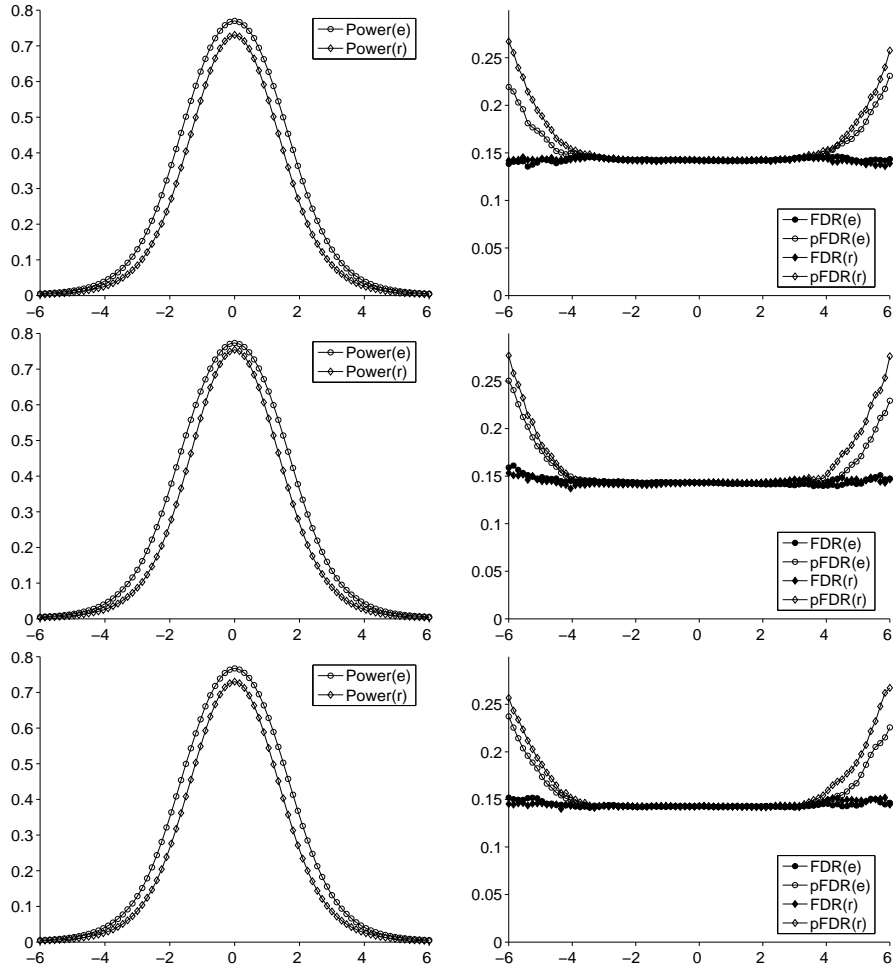

Fig 3. Power and (p)FDR *vs* log $s$ for procedures (3.9) and (4.6): group 3 (cf. Section 6.1). $a = .05$, $r = 0, 1/2, -1/2$.

## 6.2. Comparison with other procedures

In this part, we take $K \geq 2$ and compare procedures (3.9) and (4.6) with the methods in Examples 4.1, 4.4 and 4.5. The first method rejects $H_i$ with small $\prod_k \xi_{ik}$, the second one rejects $H_i$ with large $\max_k T_{ik}$ and the third one rejects $H_i$ with large $\sum_k T_{ik}$. We refer to the methods as "by-product", "by-max" and "by-sum", respectively.

The basic setup in this part is as follows. For true $H_i$, $N(\boldsymbol{\mu}_i, \Sigma_i) = N(\mathbf{0}, I_K)$ and for false $H_i$, $N(\boldsymbol{\mu}_i, \Sigma_i) = N(\boldsymbol{\mu}, \Sigma(r))$, where $\Sigma(r) = [1 + (K-1)r^2]^{-1} M'M$ with $M_{jk} = \mathbf{1}\{j = k\} + \mathbf{1}\{j \neq k\} r$. The diagonal entries of $\Sigma(r)$ are therefore 1. As in Section 6.1, FDR control parameter $\alpha = .15$, $a = .05, .02$ and $r = 0, \pm 1/5$.

TABLE 2
Parameters of the simulations in Section 6.2. The values of $\alpha_*$, $\gamma$ and $\varepsilon$ are computed for the case where the coordinates of $p$-values are independent.

|   | $\boldsymbol{\mu}$, df | $\alpha_*$ ($a = .05, .02$) | | $\boldsymbol{\gamma}$ | $\varepsilon$ |
|---|---|---|---|---|---|
| 1 | $(.5, .65, .8), 4$ | $8.58 \times 10^{-3}$ | $2.12 \times 10^{-2}$ | $(3.52, 4.93, 6.52)$ | $1/2$ |
| 2 | $(.6, .7, .8, .9, 1), 2$ | $3.25 \times 10^{-3}$ | $8.09 \times 10^{-3}$ | $(4.47, 5.48, 6.57, 7.76, 9.04)$ | $1$ |
| 3 | $(.8, .8, .8, .8), 2$ | $1.76 \times 10^{-2}$ | $4.29 \times 10^{-2}$ | $(6.57, 6.57, 6.57, 6.57)$ | $1$ |
| 4 | $(.6, .6, .6, .6), 3$ | $9.57 \times 10^{-3}$ | $2.36 \times 10^{-2}$ | $(4.27, 4.27, 4.27, 4.27)$ | $2/3$ |
| 5 | $(2, 2), 2$ | $2.88 \times 10^{-2}$ | $6.90 \times 10^{-2}$ | $(27.69, 27.69)$ | $1$ |
| 6 | $(1.5, 1.5), 3$ | $9.73 \times 10^{-3}$ | $2.40 \times 10^{-2}$ | $(16.16, 16.16)$ | $2/3$ |
| 7 | $(2, 3, 2), 10$ | $9.4 \times 10^{-19}$ | $2.4 \times 10^{-18}$ | $(39.91, 82.27, 39.91)$ | $1/5$ |

We conduct 7 groups of simulations with the values of $(\boldsymbol{\mu}, n)$ given in Table 2. Each simulation makes 3000 runs, each run tests 6000 nulls. The (p)FDR and power are computed as Monte Carlo averages of the runs.

For $K > 2$, unless $\varepsilon = 1$, the evaluation of $h(u; \boldsymbol{\nu})$ in procedure (3.9) is rather difficult. To get around this problem, by (3.10), we approximate the procedure by replacing $h(u; \boldsymbol{\nu})$ with $V_\varepsilon (u/\bar{\nu})^{K/\varepsilon}$ for all $u \in [0, 1]$. A more difficult issue is how to compare the procedures and the methods. One idea is to compare their powers at the same pFDR level $(1-a)\alpha$. However, by Examples 4.1, 4.4 and 4.5, for the values of $(\boldsymbol{\mu}, \text{df})$ in Table 2, except for the 7th one, no method attains pFDR $< \alpha = .15$. For this reason, we choose to examine the (p)FDR levels of the methods when they have the same power as either procedure at pFDR level $(1 - a)\alpha$. The steps are as follows. Take the by-product method and procedure (3.9) for example. Suppose the latter rejects $D$ false nulls when applied to $\boldsymbol{\xi}_i$. If $D > 0$, then sort $\prod_k \xi_{ik}$ in increasing order, keep rejecting the sorted nulls, starting from the first one, until $D$ false nulls are rejected; if $D = 0$, then reject no null. In this way, the number of rejected true nulls of the by-product method is minimized while the number of rejected false nulls is the same as procedure (3.9).

In each group, for each combination of $a$ and $\Sigma(r)$, procedures (3.9) and (4.6) are simulated with $\nu_k = \gamma_k$ and $c_k = (\bar{\gamma}/\gamma_k)^{1/\varepsilon}$, which are approximately the parameter values yielding maximum power for $\alpha \approx \alpha_*$. For each procedure, the by-product, by-max and by-sum methods are compared to it in the way described above. The results are reported in Tables 3–6. In groups 1–6, procedure (3.9) has more power than (4.6), often with a large margin. In all the cases where both procedures are able to control the pFDR around $(1-a)\alpha$, the methods have substantially higher FDR and pFDR when their powers are matched to that of procedure (3.9) or (4.6). The results show that the methods either cannot control the pFDR at the level of $(1 - a)\alpha$ (which is indeed the case) or, alternatively, they can only control the pFDR with much lower power than procedures (3.9) and (4.6).

Unlike groups 1–6, in group 7, each coordinate of the vector of *t*-statistics provides strong evidence to identify false nulls. By only using the 1st or 3rd coordinate, the minimum attainable pFDR is $2.4 \times 10^{-5}$ for $a = .05$ and $6 \times 10^{-5}$

Table 3

Simulation results for groups 1 and 2 described in Section 6.2. In each group, the results are organized according to $r = 0, 1/5, -1/5$. The numbers in the rows for the "by-product" method are its FDR and pFDR as its power is pegged to procedure (3.9) or (4.6). The numbers in the rows for the "by-sum" and "by-max" methods are likewise.

|  | $a = .05$ | | $a = .02$ | |
| --- | --- | --- | --- | --- |
| Group 1 | Proc. (3.9) | Proc. (4.6) | Proc. (3.9) | Proc. (4.6) |
| Power | .223 | .155 | $7.70 \times 10^{-2}$ | $5.49 \times 10^{-2}$ |
| FDR, pFDR | .141, .141 | .142, .142 | .150, .160 | .144, .162 |
| By-product | .250, .250 | .204, .204 | .250, .268 | .212, .241 |
| By-sum | .307, .307 | .283, .283 | .415, .444 | .383, .435 |
| By-max | .654, .654 | .633, .633 | .718, .768 | .660, .749 |
| Power | .351 | .316 | .210 | .181 |
| FDR, pFDR | .142, .142 | .142, .142 | .148, .148 | .147, .148 |
| By-product | .297, .297 | .264, .264 | .327, .327 | .289, .289 |
| By-sum | .333, .333 | .310, .310 | .447, .447 | .427, .428 |
| By-max | .721, .721 | .709, .709 | .836, .836 | .826, .827 |
| Power | $4.61 \times 10^{-2}$ | $3.12 \times 10^{-2}$ | $1.10 \times 10^{-2}$ | $9.28 \times 10^{-3}$ |
| FDR, pFDR | .143, .156 | .148, .165 | .141, .261 | .150, .291 |
| By-product | .194, .213 | .174, .197 | .114, .230 | .114, .250 |
| By-sum | .296, .326 | .284, .321 | .217, .440 | .204, .447 |
| By-max | .512, .564 | .485, .548 | .327, .663 | .299, .654 |
| Group 2 | Proc. (3.9) | Proc. (4.6) | Proc. (3.9) | Proc. (4.6) |
| Power | .486 | .329 | .268 | .173 |
| FDR, pFDR | .142, .142 | .143, .143 | .147, .147 | .147, .148 |
| By-product | .337, .337 | .242, .242 | .395, .396 | .313, .313 |
| By-sum | .538, .538 | .520, .520 | .734, .734 | .737, .739 |
| By-max | .809, .809 | .792, .792 | .902, .902 | .892, .894 |
| Power | .580 | .539 | .457 | .412 |
| FDR, pFDR | .142, .142 | .142, .142 | .146, .146 | .147, .147 |
| By-product | .441, .441 | .389, .389 | .510, .510 | .447, .447 |
| By-sum | .588, .588 | .567, .567 | .744, .744 | .731, .731 |
| By-max | .852, .852 | .845, .845 | .926, .926 | .923, .923 |
| Power | .316 | .138 | $3.91 \times 10^{-2}$ | $2.78 \times 10^{-2}$ |
| FDR, pFDR | .143, .143 | .142, .143 | .150, .198 | .145, .201 |
| By-product | .291, .291 | .221, .223 | .211, .287 | .189, .272 |
| By-sum | .546, .546 | .568, .572 | .558, .759 | .523, .752 |
| By-max | .785, .785 | .763, .767 | .620, .845 | .578, .830 |

for $a = .02$, and by only using the 2nd coordinate, the value is even lower. As Table 6 shows, procedures (3.9) and (4.6) identify all the false nulls. Since almost all the *p*-values of false nulls are smaller than those of true nulls, due to how the by-product method is implemented, it rejects very few true nulls and hence has near-zero (p)FDR. The same is true for the other two methods. The pFDR of procedure (3.9) is significantly lower than $(1 - a)\alpha$, because the approximation we use for $h(u; \boldsymbol{\nu})$, i.e. $V_\varepsilon(u/\bar{\nu})^{K/\varepsilon}$, is strictly greater than $h(u; \boldsymbol{\nu})$ when $u > \min \nu_k$ and hence inflates $s_i$ in (3.9). This causes the BH procedure to reject more nulls with $s_i$ not very close to 0. As these nulls are exclusively true nulls, the resulting (p)FDR is lower.

Table 4
Simulation results for groups 3 and 4 described in Section 6.2.

|  | $a = .05$ | | $a = .02$ | |
| --- | --- | --- | --- | --- |
| Group 3 | Proc. (3.9) | Proc. (4.6) | Proc. (3.9) | Proc. (4.6) |
| Power | .272 | .197 | $9.15 \times 10^{-2}$ | $6.76 \times 10^{-2}$ |
| FDR, pFDR | .142, .142 | .143, .143 | .150, .158 | .149, .162 |
| By-product | .328, .328 | .279, .279 | .344, .363 | .303, .330 |
| By-sum | .570, .570 | .572, .572 | .733, .775 | .713, .776 |
| By-max | .789, .789 | .781, .781 | .833, .880 | .802, .873 |
| Power | .444 | .409 | .291 | .257 |
| FDR, pFDR | .142, .142 | .142, .142 | .149, .149 | .150, .150 |
| By-product | .404, .404 | .368, .368 | .455, .455 | .416, .416 |
| By-sum | .584, .584 | .573, .573 | .749, .749 | .745, .745 |
| By-max | .829, .829 | .824, .824 | .913, .913 | .910, .910 |
| Power | $4.89 \times 10^{-2}$ | $3.57 \times 10^{-2}$ | $1.22 \times 10^{-2}$ | $1.01 \times 10^{-2}$ |
| FDR, pFDR | .145, .158 | .143, .160 | .151, .268 | .151, .284 |
| By-product | .234, .256 | .215, .243 | .158, .308 | .151, .321 |
| By-sum | .574, .628 | .558, .629 | .388, .757 | .357, .759 |
| By-max | .680, .744 | .651, .733 | .417, .813 | .381, .809 |
| Group 4 | Proc. (3.9) | Proc. (4.6) | Proc. (3.9) | Proc. (4.6) |
| Power | .170 | .105 | $5.69 \times 10^{-2}$ | $3.75 \times 10^{-2}$ |
| FDR, pFDR | .142, .142 | .144, .144 | .145, .161 | .144, .170 |
| By-product | .285, .285 | .227, .228 | .282, .316 | .232, .279 |
| By-sum | .440, .441 | .433, .434 | .571, .639 | .524, .629 |
| By-max | .761, .761 | .746, .747 | .763, .854 | .693, .832 |
| Power | .357 | .322 | .240 | .207 |
| FDR, pFDR | .142, .142 | .143, .143 | .147, .147 | .147, .147 |
| By-product | .363, .363 | .320, .320 | .411, .411 | .358, .358 |
| By-sum | .456, .456 | .433, .433 | .612, .612 | .594, .594 |
| By-max | .820, .820 | .812, .812 | .907, .907 | .903, .903 |
| Power | $4.26 \times 10^{-3}$ | $4.02 \times 10^{-3}$ | $2.36 \times 10^{-3}$ | $2.21 \times 10^{-3}$ |
| FDR, pFDR | .145, .294 | .149, .300 | .149, .492 | .160, .530 |
| By-product | .102, .231 | .105, .239 | $5.58 \times 10^{-2}$, .278 | $6.30 \times 10^{-2}$, .335 |
| By-sum | .217, .492 | .217, .496 | .119, .594 | .119, .635 |
| By-max | .273, .618 | .272, .621 | .143, .712 | .139, .737 |

Finally, in order to see how procedures (3.9) and (4.6) perform when the parameters $\boldsymbol{\nu}$ and $\boldsymbol{c}$ are not set to their respective asymptotically optimal values, we simulate groups 1 and 2 again, with $\boldsymbol{\nu} = \boldsymbol{c} = (1, 1, \ldots, 1)$. As Table 7 shows, across the simulations, for each procedure, the power is lower than in Table 3 but not dramatically while the (p)FDR is quite stable. The (p)FDR levels of the other 3 methods tend to be lower than in Table 3. As in group 7, this can be explained by how the methods are implemented.

TABLE 5
Simulation results for groups 5 and 6 described in Section 6.2.

|  | $a = .05$ | | $a = .02$ | |
| --- | --- | --- | --- | --- |
| Group 5 | Proc. (3.9) | Proc. (4.6) | Proc. (3.9) | Proc. (4.6) |
| Power | .772 | .733 | .372 | .342 |
| FDR, pFDR | .143, .143 | .143, .143 | .149, .149 | .148, .148 |
| By-product | .267, .267 | .246, .246 | .280, .281 | .271, .271 |
| By-sum | .387, .387 | .373, .373 | .529, .530 | .527, .528 |
| By-max | .585, .585 | .570, .570 | .697, .698 | .691, .692 |
| Power | .773 | .755 | .432 | .411 |
| FDR, pFDR | .142, .142 | .143, .143 | .149, .149 | .148, .148 |
| By-product | .283, .283 | .272, .272 | .295, .295 | .288, .288 |
| By-sum | .403, .403 | .395, .395 | .535, .536 | .532, .533 |
| By-max | .607, .607 | .599, .599 | .716, .717 | .712, .712 |
| Power | .767 | .731 | .369 | .344 |
| FDR, pFDR | .142, .142 | .142, .142 | .147, .148 | .148, .149 |
| By-product | .266, .266 | .247, .247 | .281, .282 | .272, .273 |
| By-sum | .386, .386 | .374, .374 | .529, .532 | .527, .529 |
| By-max | .584, .584 | .570, .570 | .697, .700 | .691, .695 |
| Group 6 | Proc. (3.9) | Proc. (4.6) | Proc. (3.9) | Proc. (4.6) |
| Power | .776 | .722 | .497 | .447 |
| FDR, pFDR | .143, .143 | .143, .143 | .146, .146 | .147, .147 |
| By-product | .239, .239 | .207, .207 | .259, .259 | .237, .237 |
| By-sum | .295, .295 | .269, .269 | .391, .391 | .377, .377 |
| By-max | .540, .540 | .513, .513 | .654, .654 | .639, .639 |
| Power | .771 | .747 | .538 | .510 |
| FDR, pFDR | .142, .142 | .142, .142 | .148, .148 | .149, .149 |
| By-product | .258, .258 | .240, .240 | .278, .278 | .263, .263 |
| By-sum | .315, .315 | .301, .301 | .410, .410 | .399, .399 |
| By-max | .571, .571 | .557, .557 | .685, .685 | .676, .676 |
| Power | .773 | .711 | .476 | .425 |
| FDR, pFDR | .143, .143 | .143, .143 | .148, .148 | .148, .148 |
| By-product | .233, .233 | .200, .200 | .257, .257 | .235, .235 |
| By-sum | .288, .288 | .262, .262 | .388, .388 | .375, .375 |
| By-max | .530, .530 | .501, .501 | .646, .646 | .631, .631 |

## 7. Discussion

### 7.1. Role of p-values

We have followed the tradition of using *p*-values for hypothesis testing. The general procedure in the work, i.e., (3.2), utilizes the fact that the *p*-value of a continuous multivariate statistic can be defined in such a way that its coordinates are i.i.d. $\sim \text{Unif}(0,1)$. The interpretation of *p*-value as a measure on how "rare" or "suspi0cious" an observation looks is irrelevant, even though in many cases smaller *p*-values are indeed more likely to be associated with false nulls.

TABLE 6
Simulation results for group 7 described in Section 6.2. "–" means value equal to the nonmissing value in the same row.

| Group 7 | $a = .05$ | | $a = .02$ | |
| --- | --- | --- | --- | --- |
| | Proc. (3.9) | Proc. (4.6) | Proc. (3.9) | Proc. (4.6) |
| Power | 1 | 1 | 1 | 1 |
| FDR, pFDR | .103, .103 | .142, .142 | .118, .118 | .146, .146 |
| By-product | 0, – | –, – | 0, – | –, – |
| By-sum | $1.2 \times 10^{-6}$, – | –, – | $5.0 \times 10^{-6}$, – | –, – |
| By-max | $3.1 \times 10^{-3}$, – | –, – | $5.4 \times 10^{-3}$, – | –, – |
| Power | 1 | 1 | | |
| FDR, pFDR | .103, .103 | .143, .143 | .118, .118 | .148, .148 |
| By-product | $1.2 \times 10^{-6}$, – | –, – | 0, – | –, – |
| By-sum | $3.6 \times 10^{-6}$, – | –, – | $1.3 \times 10^{-5}$, – | –, – |
| By-max | $5.7 \times 10^{-3}$, – | –, – | $9.6 \times 10^{-3}$, – | –, – |
| Power | 1 | 1 | 1 | 1 |
| FDR, pFDR | .103, .103 | .143, .143 | .119, .119 | .148, .148 |
| By-product | 0, – | –, – | 0, – | –, – |
| By-sum | 0, – | –, – | $2.7 \times 10^{-6}$, – | –, – |
| By-max | $2.9 \times 10^{-3}$, – | –, – | $4.9 \times 10^{-3}$, – | –, – |

Thus, in this work, *p*-values serve as a mechanism to "flatten" the probability landscape of true nulls and hence facilitates exploring subtle differences between true and false nulls.

Since what essentially matters to procedure (3.2) is nested *events* with specific probabilities, it can be easily modified to directly handle test statistics instead of their *p*-values. Indeed, in (3.2), $D_t \in [0,1]^K$ can be replaced with nested $E_t$ in the domain of the test statistics, such that $P(E_t) = t$ under true nulls. Analysis on the power of the modification might yield some useful insight. For example, weighted $L^2$ norms are commonly used as criterion for acceptance/rejection. However, as shown in Examples 4.5 and 4.6, in more challenging cases, one may need to consider $L^p$ norms with $p < 0$. On the other hand, the modification does not simplify the testing problem, as probabilities still have to be evaluated. Nevertheless, as remarked next, the notion of using nested regions in spaces other than $[0,1]^K$ is useful.

### 7.2. Incorporating discrete components

Often times, test statistics have nontrivial discrete components. For example, test statistics for different nulls may have different dimensions or degrees of freedom. In this case, the discrete component may be expressed as a scalar. However, if the test statistics are multivariate but only partially observed, then the discrete component in general have to be set-valued accounting for observed coordinates. Procedure (3.2) can be modified as follows. Suppose $Z$ is the discrete component of test statistic $\boldsymbol{T}$ such that for any $z$, the conditional distribution of $\boldsymbol{T}$ given $Z = z$ has a density and is $K(z)$ dimensional. Then, in

TABLE 7

Simulation results for groups 1 and 2. The setting is similar to that in Table 3, except that $\boldsymbol{\nu}$ in (3.9) and $\boldsymbol{c}$ in (4.6) are set equal to $(1, 1, \ldots, 1)$ instead of according to $\boldsymbol{\gamma}$.

|  | $a = .05$ | | $a = .02$ | |
| --- | --- | --- | --- | --- |
| Group 1 | Proc. (3.9) | Proc. (4.6) | Proc. (3.9) | Proc. (4.6) |
| Power | .206 | .138 | $6.94 \times 10^{-2}$ | $4.88 \times 10^{-2}$ |
| FDR, pFDR | .143, .143 | .143, .143 | .146, .159 | .144, .163 |
| By-product | .240, .240 | .191, .191 | .236, .257 | .203, .231 |
| By-sum | .302, .302 | .278, .278 | .405, .441 | .378, .430 |
| By-max | .650, .650 | .626, .627 | .700, .762 | .657, .746 |
| Power | .340 | .293 | .199 | .166 |
| FDR, pFDR | .143, .143 | .143, .143 | .147, .147 | .146, .146 |
| By-product | .286, .286 | .242, .242 | .315, .315 | .271, .272 |
| By-sum | .326, .326 | .296, .296 | .443, .443 | .419, .420 |
| By-max | .717, .717 | .701, .701 | .835, .835 | .823, .824 |
| Power | $3.78 \times 10^{-2}$ | $2.55 \times 10^{-2}$ | $1.06 \times 10^{-2}$ | $8.55 \times 10^{-3}$ |
| FDR, pFDR | .140, .155 | .142, .166 | .146, .279 | .144, .289 |
| By-product | .175, .196 | .154, .182 | .107, .227 | .108, .247 |
| By-sum | .284, .318 | .267, .316 | .210, .446 | .196, .448 |
| By-max | .495, .554 | .453, .537 | .313, .666 | .286, .654 |
| Group 2 | Proc. (3.9) | Proc. (4.6) | Proc. (3.9) | Proc. (4.6) |
| Power | .458 | .300 | .248 | .155 |
| FDR, pFDR | .142, .142 | .142, .142 | .147, .147 | .148, .148 |
| By-product | .317, .317 | .226, .226 | .378, .378 | .296, .297 |
| By-sum | .533, .533 | .519, .519 | .734, .734 | .739, .741 |
| By-max | .806, .806 | .788, .788 | .901, .901 | .889, .892 |
| Power | .562 | .508 | .441 | .388 |
| FDR, pFDR | .142, .142 | .142, .142 | .145, .145 | .145, .145 |
| By-product | .420, .420 | .353, .353 | .483, .483 | .412, .412 |
| By-sum | .581, .581 | .554, .554 | .738, .738 | .723, .723 |
| By-max | .850, .850 | .840, .840 | .925, .925 | .920, .920 |
| Power | .251 | .112 | $3.34 \times 10^{-2}$ | $2.35 \times 10^{-2}$ |
| FDR, pFDR | .141, .141 | .141, .142 | .149, .204 | .146, .212 |
| By-product | .266, .266 | .207, .209 | .200, .281 | .176, .264 |
| By-sum | .551, .552 | .574, .580 | .537, .755 | .500, .754 |
| By-max | .778, .780 | .757, .765 | .595, .838 | .550, .829 |

(3.2), redefine $D_t$ as a nested subsets in the disjoint union of $[0,1]^{K(z)}$, such that $\sum_z \ell(D_t \cap [0,1]^{K(z)})p_z = t$, where $p_z$ is the probability of $Z = z$ under true nulls. The analysis in previous sections still works and requires no substantial extra changes.

An apparently simpler alternative is to conduct separate tests on statistics with different values of the discrete components. This alternative fails to take into account the distribution of the discrete components and hence may have lower power.

### 7.3. Power optimization

When the distribution under false nulls is only partially known, it can be a difficult issue how to attain maximum power. To see this better, consider testing the null "$\boldsymbol{\mu} = 0$" for $N(\boldsymbol{\mu}, I)$ based on a single observation $\boldsymbol{X}$. As the variance is *known* to be $I$, the most powerful test statistic would be $\boldsymbol{\nu}'\boldsymbol{X}$, provided that the true value $\boldsymbol{\nu}$ of $\boldsymbol{\mu}$ under false nulls is known. However, when $\boldsymbol{\nu}$ is unknown, unless there is strong evidence on its whereabouts, one has to search in a large region of $\boldsymbol{\mu}$ to improve the power, which becomes more difficult as the dimension of $\boldsymbol{\nu}$ gets higher.

One way to improve power is to restrict the search to parametric families of nested regions. This is the approach taken in Section 3.3. If the parameter involved is of high dimension, some type of stochastic optimization [19] may be needed. On the other hand, regions that attain maximum power may consist of several disconnected regions, which makes it difficult to use a single parametric family of nested regions to approximate them. An alternative way therefore is to try different families of nested regions at different locations in the domain of $p$-values and combine the results appropriately [7].

### Appendix

In this section, we shall denote $I = [0, 1]^K$.

### A.1. Theoretical details for Section 3

*Proof of Proposition 3.1.* Since $D_0 = \emptyset$ and $D_1 = I$, it suffices to show that $D_t$, $t \in (0, 1)$ satisfy (3.1). Observe that $h^{-1}$ is continuous and strictly increasing on $(0, 1)$. Then, as $\Gamma_u$ is right-continuous, $D_t$ is right-continuous. It is clear that $D_t$ is increasing and $\ell(D_t) = \ell(\Gamma_{h^*(t)}) = h(h^*(t)) = t$. □

*Proof of Lemma 3.1.* 1) The following "sandwiched convergence" is needed: if $0 \leq a_n(x) \leq b_n(x)$ such that $a_n(x) \to a(x)$, $b_n(x) \to b(x)$ a.e. and $\int b_n \to \int b < \infty$, then $\int a_n \to \int a$. For each $k$, denote by $g(x_1, \ldots, x_k)$

$$\int_{-\infty}^{x_k} f_k(x_1, \ldots, x_{k-1}, z)\, dz = \int f_k(x_1, \ldots, x_{k-1}, z)\mathbf{1}\{z \leq x_k\}\, dz.$$

The function in the second integral is dominated by $f_k(x_1, \ldots, x_{k-1}, z)$. If $(x_1, \ldots, x_k) \to (y_1, \ldots, y_k)$, then, by the continuity of $f_k$ and $f_{k-1}$,

$$f_k(x_1, \ldots, x_{k-1}, z)\mathbf{1}\{z \leq x_k\} \to f_k(y_1, \ldots, y_{k-1}, z)\mathbf{1}\{z \leq y_k\}, \text{ for } z \neq y_k$$
$$f_k(x_1, \ldots, x_{k-1}, z) \to f_k(y_1, \ldots, y_{k-1}, z),$$
$$\int f_k(x_1, \ldots, x_{k-1}, z)\, dz = f_{k-1}(x_1, \ldots, x_{k-1}) \to \int f_k(y_1, \ldots, y_{k-1}, z)\, dz.$$

By the sandwiched convergence, $g(x_1, \ldots, x_k) \to g(y_1, \ldots, y_k)$. Thus, $g$ is continuous. As $\phi_k(\boldsymbol{x}) = g(x_1, \ldots, x_k)/f_{k-1}(x_1, \ldots, x_{k-1})$ for $\boldsymbol{x} \in \text{sppt}(q_0)$, $\phi_k(\boldsymbol{x})$ is continuous in $\text{sppt}(q_0)$. Therefore, $\phi \in C(\text{sppt}(q_0))$.

Let $\boldsymbol{x}, \boldsymbol{y} \in \text{sppt}(q_0)$. Suppose $\boldsymbol{x} \ne \boldsymbol{y}$ such that $x_i = y_i$ for $i < k$ and $x_k < y_k$. Then $\phi_k(\boldsymbol{x}) \le \phi_k(\boldsymbol{y})$. Assume equality holds. Then

$$0 = [\phi_k(\boldsymbol{y}) - \phi_k(\boldsymbol{x})] f_{k-1}(x_1, \ldots, x_{k-1})$$
$$= \int \cdots \int \int_{x_k}^{y_k} q_0(x_1, \ldots, x_{k-1}, z, u_{k+1}, \ldots, u_K) \, dz \, du_{k+1} \cdots du_K.$$

Since $q_0$ is continuous, the above formula implies $q_0(x_1, \ldots, x_{k-1}, z, u_{k+1}, \ldots, u_K) = 0$ for $z \in [x_k, y_k]$ and $u_{k+1}, \ldots, u_K \in \mathbb{R}$, in particular, $q_0(\boldsymbol{x}) = q_0(\boldsymbol{y}) = 0$. The contradiction implies $\phi_k(\boldsymbol{x}) < \phi_k(\boldsymbol{y})$ and so $\phi(\boldsymbol{x}) \ne \phi(\boldsymbol{y})$.

2) Let $\boldsymbol{X} \sim Q_0$. Then $P(\boldsymbol{X} \in \text{sppt}(q_0)) = 1$. Since $\phi \in C(\text{sppt}(q_0))$, $\boldsymbol{\xi} = \phi(\boldsymbol{X})$ is a well-defined random variable. For $\boldsymbol{x} \in \text{sppt}(q_0)$, by 1), conditional on $X_i = x_i$, $i < k$, $X_k$ has a continuous distribution and hence $\xi_k \sim \text{Unif}(0, 1)$. Since the conditional distribution of $\xi_k$ is the same regardless of $x_1, \ldots, x_{k-1}$, $\xi_k$ is independent of $X_1, \ldots, X_{k-1}$ and thus independent of $\xi_1, \ldots, \xi_{k-1}$. This gives $\boldsymbol{\xi} \sim \text{Unif}(I)$.

3) Let $\boldsymbol{X} \sim Q_1$. As in 2), $\boldsymbol{\xi} = \phi(\boldsymbol{X})$ is a well-defined random variable. Denote $r = q_1/q_0$. For $\boldsymbol{t} \in \mathbb{R}^K$,

$$E[e^{i\boldsymbol{t}'\boldsymbol{\xi}}] = E_{Q_0}\left[e^{i\boldsymbol{t}'\boldsymbol{\xi}} r(\boldsymbol{X})\right] = E_{Q_0}\left[e^{i\boldsymbol{t}'\boldsymbol{\xi}} r(\phi^{-1}(\boldsymbol{\xi}))\right],$$

where the first equality holds since $\text{sppt}(q_1) \subset \text{sppt}(q_0)$ and $r(\boldsymbol{X})$ is a well-defined random variable due to $r \in C(\text{sppt}(q_0))$. Since $\phi$ is 1-to-1 and continuous on $\text{sppt}(q_0)$, $E := \phi(\text{sppt}(q_0))$ is open and $\phi^{-1} \in C(E)$ [15]. As $\ell(I \setminus E) = Q_0(\boldsymbol{\xi} \notin E) = 0$, $r(\phi^{-1}(\boldsymbol{x}))$ is Borel measurable on $I$ and by 2), the last expectation equals $\int_I e^{i\boldsymbol{t}'\boldsymbol{u}} r(\phi^{-1}(\boldsymbol{u})) \, d\boldsymbol{u}$. Thus, the characteristic function of $\boldsymbol{\xi}$ under $Q_1$ is the same as that of a random variable with density $r(\phi^{-1}(\boldsymbol{u}))$, $\boldsymbol{u} \in I$. Since $r(\phi^{-1}(\boldsymbol{u})) \in C(E)$, this proves 3). $\square$

To show Proposition 3.2, we need a few preliminary results. Recall $p_i = h(g(\boldsymbol{\xi}_i))$.

**Lemma A.1.1.** *Let $h$ be continuous. Then 1) $p_i \le t \iff \boldsymbol{\xi}_i \in D_t$ and 2) under the assumptions of Lemma 3.1,*

$$P(p_i \le t) = \begin{cases} t & \text{if } H_i \text{ is true} \\ G(D_t) & \text{if } H_i \text{ is false} \end{cases}, \quad t \in [0, 1],$$

*and $G(D_t)$ is strictly concave.*

*Proof.* It can be seen that $p_i$ is a well-defined random variable and for $t \in [0, 1]$, $h(h^*(t)) = t$. Then $p_i \le t \iff g(\boldsymbol{\xi}_i) \ge h^*(t) \iff \boldsymbol{\xi}_i \in D_t$. Under true $H_i$, $P(p_i \le t) = \ell(D_t) = t$. Under false $H_i$, $P(p_i \le t) = \int_{D_t} g = G(D_t)$. Given $0 \leq t_1 < t_2 < t_3 \leq 1$, let $u_k = h^*(t_k)$. By the continuity of $h$, $t_k = h(u_k)$ and $u_1 > u_2 > u_3$. As $D_{t_k} = \Gamma_{u_k}$,

$$r_k := \frac{G(D_{t_{k+1}}) - G(D_{t_k})}{t_{k+1} - t_k} = \frac{1}{h(u_{k+1}) - h(u_k)} \int_{\Gamma_{u_{k+1}} \setminus \Gamma_{u_k}} g(\boldsymbol{x}) \, d\boldsymbol{x}.$$

Since $h(u_k) = \ell(\Gamma_{u_k})$ and $g(\boldsymbol{x}) \in (u_{k+1}, u_k]$ on $\Gamma_{u_{k+1}} \setminus \Gamma_{u_k}$, $u_{k+1} < r_k \leq u_k$. As a result, $r_1 > r_2$. Therefore, the distribution of $p_i$ is strictly concave. □

**Lemma A.1.2.** *Let $\eta_1, \ldots, \eta_n$ be independent Bernoulli random variables such that $p_i = P(\eta_i = 1)$ are decreasing. Let $S = \eta_1 + \cdots + \eta_n$. Then $E[\eta_i/(S \vee 1)]$ is decreasing.*

*Proof.* Let $i < j$. Then $\eta_i, \eta_j$ and $X = S - \eta_i - \eta_j$ are independent, giving

$$E\left[\frac{\eta_i}{S \vee 1}\right] = p_i E\left[\frac{1}{1 + \eta_j + X}\right], \quad E\left[\frac{\eta_j}{S \vee 1}\right] = p_j E\left[\frac{1}{1 + \eta_i + X}\right].$$

Since $p_i \geq p_j$, $(1+\eta_j+X)^{-1}$ stochastically dominates $(1+\eta_i+X)^{-1}$. Therefore, $p_i E[(1 + \eta_j + X)^{-1}] \geq p_j E[(1 + \eta_i + X)^{-1}]$. □

**Lemma A.1.3.** *Let $s_n, \eta_1, \ldots, \eta_n \in [0,1]$ be jointly distributed, such that $s_n \xrightarrow{P} s \in [0,1]$ as $n \to \infty$ and $\eta_i$ are i.i.d. $\sim F$. Let $F_n$ be the empirical distribution of $\eta_1, \ldots, \eta_n$. If $F$ is continuous and strictly increasing on $[0,1]$, then $F_n^*(s_n) \xrightarrow{P} F^*(s)$.*

*Proof.* Recall $\sup|F_n - F| \xrightarrow{P} 0$. Let $x_n = F_n^*(s_n)$ and $x = F^*(s)$. Since $F$ is continuous, $s = F(x)$. Suppose $x \in (0,1)$. Given $\epsilon \in (0,x)$, by $F_n(x_n - \epsilon) < s_n$, $\{x_n > x + 2\epsilon\} \subset \{s_n > F_n(x+\epsilon)\}$. By $s_n \xrightarrow{P} s = F(x)$ and $F_n(x+\epsilon) \xrightarrow{P} F(x+\epsilon) > F(x)$, $P(\{x_n > x + 2\epsilon\}) \to 0$. On the other hand, $\{x_n < x - \epsilon\} \subset \{s_n \leq F_n(x-\epsilon)\}$. By $F_n(x-\epsilon) \xrightarrow{P} F(x-\epsilon) < F(x)$, $P(\{x_n < x - \epsilon\}) \to 0$. Therefore, $x_n \xrightarrow{P} 0$. The case where $x = 0$ or $1$ is similarly proved. □

*Proof of Proposition 3.2.* First, since $h \in C$ is decreasing, $\cap_{s>t} D_s = \cap_{s>t} \Gamma_{h^*(s)} = \cap_{s>t} \{g(\boldsymbol{x}) \geq h^*(s)\} = \cap_{s>t} \{h(g(\boldsymbol{x})) \leq s\} = \{h(g(\boldsymbol{x})) \leq t\} = \{g(\boldsymbol{x}) \geq h^*(t)\} = D_t$, proving the right-continuity of $D_t$. By the continuity of $h$, $\ell(D_t) = \ell(\Gamma_{h^*(t)}) = h(h^*(t)) = t$. The rest of (3.1) is easy to check.

Denote by $N$ the number of true nulls, and for any given procedure, denote by $R$ and $V$ the numbers of rejected nulls and rejected true nulls, respectively. We shall show i) procedure (3.2) with $D_t$ satisfies conditions (A) and (B) and attains FDR $= (1-a)\alpha$; ii) the search for procedures with maximum power can be restricted to those that reject and only reject nulls with largest $g(\boldsymbol{\xi}_i)$; and iii) for such a procedure, $R/n$ converges in probability to a nonrandom number. From these results, the proof will follow without much difficulty.

i) By Lemma A.1.1 1), procedure (3.2) with $D_t = \Gamma_{h^*(t)}$ is the same as the BH procedure applied to $p_1, \ldots, p_n$. Therefore, statement 1) holds and FDR $= (1 - a)\alpha$ [2, 22]. Since the set of rejected nulls is uniquely determined by $\boldsymbol{\xi}_1, \ldots, \boldsymbol{\xi}_n$, the procedure satisfies condition (A).

Recall $\alpha_* = 1/(1 - a + a \sup g)$. By Lemma A.1.1 2), $F(t) = (1-a)t + aG(D_t)$ is strictly concave. Then

$$\alpha_* = \lim_{u \to \sup g} \frac{1}{1 - a + aG(\Gamma_u)/\ell(\Gamma_u)} = \frac{1}{F'(0)}.$$

For $\alpha \in (\alpha_*, 1)$, $t/\alpha = F(t)$ has a unique positive solution $t_* \in (0, 1)$. Since $p_i$ are i.i.d. $\sim F$, by [12], for $\tau$ in (3.2), $\tau \xrightarrow{P} t_*$ as $n \to \infty$ and procedure (3.2) asymptotically has the same power as the one that rejects $H_i$ with $p_i \le t_*$. By the law of large numbers, $R/n \xrightarrow{P} F(t_*)$. On the other hand, for $\alpha < \alpha_*$, by [7], $R/n \xrightarrow{P} 0$. In either case, condition (B) is satisfied.

ii) Given $\boldsymbol{\xi}_1, \ldots, \boldsymbol{\xi}_n$, $\theta_1, \ldots, \theta_n$ are independent Bernoulli variables with

$$r_i := P(\theta_i = 1 \mid \boldsymbol{\xi}_1, \ldots, \boldsymbol{\xi}_n) = \frac{ag(\boldsymbol{\xi}_i)}{1 - a + ag(\boldsymbol{\xi}_i)}.$$

Sort $r_i$ into $r_{(1)} \ge r_{(2)} \ge \ldots \ge r_{(n)}$. Let $\boldsymbol{\delta}$ be a procedure satisfying condition (A). Then $R = \sum_{i=1}^n \delta_i$ and $V = \sum_{i=1}^n (1 - \theta_i)\delta_i$. By condition (A), $\theta_i$ is conditionally independent of $(\delta_i, R)$ given $\boldsymbol{\xi}_1, \ldots, \boldsymbol{\xi}_n$. Then

$$E[V \mid \boldsymbol{\xi}_1, \ldots, \boldsymbol{\xi}_n, R] = \sum_{i=1}^n E[(1 - \theta_i)\delta_i \mid \boldsymbol{\xi}_1, \ldots, \boldsymbol{\xi}_n, R]$$

$$= \sum_{i=1}^n (1 - r_i) E[\delta_i \mid \boldsymbol{\xi}_1, \ldots, \boldsymbol{\xi}_n, R] \ge \sum_{i=1}^R (1 - r_{(i)}),$$

where the last inequality is due to $R = \sum_{i=1}^n \delta_i$. Then

$$\text{FDR} = E\left[\frac{E[V \mid \boldsymbol{\xi}_1, \ldots, \boldsymbol{\xi}_n, R]}{R \vee 1}\right] \ge E\left[\frac{1}{R \vee 1} \sum_{i=1}^R (1 - r_{(i)})\right]$$

with equality being true if rejected nulls are exactly those with the $R$ largest $r_i$. On the other hand, since $N = \sum_{i=1}^n \theta_i$, by Lemma A.1.2,

$$E\left[\frac{R - V}{N \vee 1} \mid \boldsymbol{\xi}_1, \ldots, \boldsymbol{\xi}_n, R\right] = \sum_{i=1}^n E\left[\frac{\theta_i \delta_i}{N \vee 1} \mid \boldsymbol{\xi}_1, \ldots, \boldsymbol{\xi}_n, R\right]$$

$$= \sum_{i=1}^n E\left[\frac{\theta_i}{N \vee 1} \mid \boldsymbol{\xi}_1, \ldots, \boldsymbol{\xi}_n\right] E[\delta_i \mid \boldsymbol{\xi}_1, \ldots, \boldsymbol{\xi}_n, R]$$

$$\le \sum_{i=1}^R E\left[\frac{\theta_{(i)}}{N \vee 1} \mid \boldsymbol{\xi}_1, \ldots, \boldsymbol{\xi}_n\right]$$

where $\theta_{(i)}$ corresponds to the null with the $i$th largest $r_i$. Then

$$\text{power} \leq E\left[\sum_{i=1}^{R} E\left[\frac{\theta_{(i)}}{N \vee 1} \,\bigg|\, \boldsymbol{\xi}_1, \ldots, \boldsymbol{\xi}_n\right]\right].$$

Note that $r_i \geq r_j \iff g(\boldsymbol{\xi}_i) \geq g(\boldsymbol{\xi}_j) \iff p_i \leq p_j$. Construct procedure $\boldsymbol{\delta}'$ which first applies $\boldsymbol{\delta}$ and then, provided $\boldsymbol{\delta}$ rejects $R$ nulls, rejects nulls with the $R$ smallest $p_i$ instead. It follows that 1) if $\boldsymbol{\delta}$ has FDR $\leq (1-a)\alpha$, then so does $\boldsymbol{\delta}'$; 2) $\boldsymbol{\delta}'$ is at least as powerful as $\boldsymbol{\delta}$; 3) if $\boldsymbol{\delta}$ satisfies condition (A), then, as the second step of $\boldsymbol{\delta}'$ is uniquely determined by $\boldsymbol{\xi}_1, \ldots, \boldsymbol{\xi}_n$, $\boldsymbol{\delta}'$ satisfies condition (A) as well; and 4) since $\boldsymbol{\delta}$ and $\boldsymbol{\delta}'$ reject the same number of nulls, if one satisfied condition (B), the other does as well.

iii) Let $\boldsymbol{\delta}$ satisfy conditions (A), (B) and attain maximum power asymptotically while having FDR $\leq (1-a)\alpha$. As $n \to \infty$, the empirical distribution of $\boldsymbol{\xi}_1, \ldots, \boldsymbol{\xi}_n$ converges to the distribution that has density $1 - a + ag(\boldsymbol{x})$. By condition (B), $R/n \xrightarrow{\text{P}}$ some $t_* \in [0,1]$. Let $F_n$ be the empirical distribution of $p_i$. Then $F_n \xrightarrow{\text{P}} F$. Since $F$ is strictly concave and continuous, $F$ is strictly increasing. By Lemma A.1.3, $F_n^*(R/n) \xrightarrow{\text{P}} F^*(s)$. Since $\boldsymbol{\delta}$ rejects and only rejects nulls with $p_i \leq F_n^*(R/n)$, $\boldsymbol{\delta}$ is asymptotically equivalent to a procedure which rejects and only rejects nulls with $p_i \leq t_* = F^*(s)$. If $t_* > 0$, by the law of large numbers and dominate convergence, as $n \to \infty$,

$$\text{FDR} = (1 + o(1))E\left[\frac{\#\{i : p_i \leq t_*, \theta_i = 0\}}{\#\{i : p_i \leq t_*\} \vee 1}\right] \to \frac{(1-a)t_*}{F(t_*)},$$

$$\text{power} = (1 + o(1))E\left[\frac{\#\{i : p_i \leq t_*, \theta_i = 1\}}{\#\{i : \theta_i \leq t_*\} \vee 1}\right] \to G(D_{t_*}).$$

In order to attain maximum power while maintaining FDR $\leq (1-a)\alpha$, $t_*$ has to be the largest value of $t$ satisfying $t/F(t) \leq \alpha$. It is easy to see that for $\alpha \in (\alpha_*, 1)$, $t_*$ is the unique positive solution of $t/\alpha = F(t)$. Combined with part i) of the proof, this shows procedure (3.2) with $D_t = \Gamma_{h^*(t)}$ can be taken as $\boldsymbol{\delta}$. Furthermore, in this case, power $\to G(D_{t_*}) > 0$ and since $P(R > 0) \to 1$, pFDR $= (1 + o(1))\text{FDR} \to (1-a)\alpha$. Thus 2) is proved.

On the other hand, for $\alpha < \alpha_*$, no $t > 0$ satisfies $t/F(t) \leq \alpha$. As a result, $t_* = 0$. Thus, the power of $\boldsymbol{\delta}$ is asymptotically 0 and procedure (3.2) with $D_t$ again can be taken as $\boldsymbol{\delta}$. Furthermore, by [7], the procedure has pFDR $\to (1-a)\alpha_*$. This proves 3). $\square$

*Proof of Lemma 3.2.* By change of variable $x_k = (u/\nu_k)^{1/\varepsilon} z_k$, for $0 < u \leq \min \nu_k$,

$$h(u; \boldsymbol{\nu}) = \int_I \mathbf{1}\left\{\boldsymbol{\nu}' \boldsymbol{x}^\varepsilon \leq u\right\} d\boldsymbol{x} = \int_I \mathbf{1}\left\{\sum z_k^\varepsilon \leq 1\right\} \prod \left(\frac{u}{\nu_k}\right)^{1/\varepsilon} d\boldsymbol{z},$$

which yields (3.10). Likewise,

$$\int_I x_k^\varepsilon \mathbf{1}\{\boldsymbol{\nu}'\boldsymbol{x}^\varepsilon \leq u\}\, d\boldsymbol{x} = \frac{\bar{\nu}}{\nu_k}\left(\frac{u}{\bar{\nu}}\right)^{K/\varepsilon+1}\int_I z_k^\varepsilon \mathbf{1}\left\{\sum z_k^\varepsilon \leq 1\right\} d\boldsymbol{z}.$$

By symmetry,

$$\int_I z_k^\varepsilon \mathbf{1}\left\{\sum z_k^\varepsilon \leq 1\right\} d\boldsymbol{z} = \frac{1}{K}\int_I \left(\sum z_k^\varepsilon\right)\mathbf{1}\left\{\sum z_k^\varepsilon \leq 1\right\} d\boldsymbol{z}.$$

By $V_\varepsilon(u) := \ell(\{\boldsymbol{z} \in \mathbb{R}^K : z_k \geq 0, \sum z_k^\varepsilon \leq u\}) = V_\varepsilon u^{K/\varepsilon}$ and change of measure, the right hand side equals $(1/K)\int_0^1 t\, dV_\varepsilon(t) = V_\varepsilon/(K+\varepsilon)$. It is then easy to get (3.11). Finally, by (3.7),

$$\int_{\Gamma_u(\boldsymbol{\nu})} g = g(\mathbf{0})\left(h(u;\boldsymbol{\nu}) - \sum \gamma_k \int_{\Gamma_u(\boldsymbol{\nu})} x_k^\varepsilon\, d\boldsymbol{x}\right) + g(\mathbf{0})\int_{\Gamma_u(\boldsymbol{\nu})} r.$$

As $u \to 0$, $\sup_{\Gamma_u(\boldsymbol{\nu})}|\boldsymbol{x}| \to 0$, implying $r(\boldsymbol{x})/|\boldsymbol{x}|^\varepsilon \to 0$ uniformly on $\Gamma_u(\boldsymbol{\nu})$, and by (3.11), $\int_{\Gamma_u(\boldsymbol{\nu})} r = o(u^{K/\varepsilon+1})$, which together with (3.10) yields (3.12). □

*Proof of (3.13).* Let $\nu_1 < \nu_2$ without loss of generality. For $0 < u < \nu_1 + \nu_2$, $h(u;\boldsymbol{\nu})$ is equal to

$$\int_0^1 \mathbf{1}\{\nu_1 x^\varepsilon \leq u - \nu_2\}\, dx + \int_0^1 \mathbf{1}\{u - \nu_2 < \nu_1 x^\varepsilon \leq u\}\left(\frac{u - \nu_1 x^\varepsilon}{\nu_2}\right)^{1/\varepsilon} dx.$$

The first integral on the right hand side equals $[(u-\nu_2)/\nu_1]^{1/\varepsilon}$ if $u \geq \nu_2$ and 0 otherwise. By variable substitution $z = \nu_1 x^\varepsilon/u$, the second integral equals

$$\frac{1}{\varepsilon}\left(\frac{u^2}{\nu_1\nu_2}\right)^{1/\varepsilon}\int_0^{\nu_1/u} \mathbf{1}\left\{1 - \frac{\nu_2}{u} < z \leq 1\right\} z^{1/\varepsilon-1}(1-z)^{1/\varepsilon}\, dz$$

$$= \frac{\Gamma(1/\varepsilon)^2}{2\varepsilon\Gamma(2/\varepsilon)}\left(\frac{u^2}{\nu_1\nu_2}\right)^{1/\varepsilon}\left[F\left(\frac{\nu_1}{u} \wedge 1\right) - F\left(\left(1 - \frac{\nu_2}{u}\right) \vee 0\right)\right],$$

where $F$ is the Beta distribution function with parameters $1/\varepsilon$ and $1/\varepsilon + 1$. Since $F(x \wedge 1) = F(x)$ and $F(x \vee 0) = F(x)$, this yields the proof. □

### A.2. Theoretical details for Section 4

*Proof of Proposition 4.1.* Denote $h(u) = h(u;\boldsymbol{\nu})$ and $D_t$ the regularization of $\Gamma_u(\boldsymbol{\nu})$ defined in (3.8). Let $Z_i = \boldsymbol{\nu}'\boldsymbol{\xi}_i^\varepsilon$. Under true $H_i$, since $\boldsymbol{\xi}_i \sim \text{Unif}(I)$, $P(Z_i \leq u) = h(u)$. Since $h$ is continuous, $h(Z_i) \sim \text{Unif}(0,1)$. On the other hand, under false $H_i$, $P(h(Z_i) \leq t) = P(\boldsymbol{\nu}'\boldsymbol{\xi}_i^\varepsilon \leq h^*(t)) = G(\Gamma_{h^*(t)}(\boldsymbol{\nu})) = G(D_t)$, so by (4.2), the density of $h(Z_i)$ at 0 is $g(\mathbf{0})$. Because procedure (3.9) is the BH procedure applied to $h(Z_i)$, it follows that its minimum attainable pFDR is $(1-a)\alpha_*$.

To get (4.1), let $F(t) = (1-a)t + aG(D_t)$ and $t_*$ the maximum solution to $F(t) = t/\alpha$. Then

$$G(D_{t_*}) = \left(\frac{1}{\alpha a} - \frac{1}{a} + 1\right) t_* = \left(\frac{1}{\alpha a} - \frac{1}{\alpha_* a} + g(\mathbf{0})\right) t_*.$$

Replacing the left hand side by (4.2), it is seen that as $\alpha \downarrow \alpha_*$,

$$\frac{(t/V_\varepsilon)^{\varepsilon/K} \bar{\nu}}{K+\varepsilon} \sum \frac{\gamma_k}{\nu_k} \sim \frac{\alpha - \alpha_*}{a\alpha_*^2 g(\mathbf{0})}$$

$$\implies t_* \sim V_\varepsilon \left(\frac{K+\varepsilon}{a\alpha_*^2 g(\mathbf{0})} \bigg/ \sum \frac{\bar{\nu}\gamma_k}{\nu_k}\right)^{K/\varepsilon} (\alpha - \alpha_*)^{K/\varepsilon}.$$

Clearly, as $\alpha \downarrow \alpha_*$, $t_* \to 0$. If $G(D_t)$ is strictly concave, then by [12], $\text{Pow}(\alpha) = G(D_{t_*}) \sim g(\mathbf{0}) t_*$, which, combined with the asymptotic of $t_*$, proves (4.1). However, it is not clear whether $G(D_t)$ is strictly concave in general. To get around the problem, we use the following argument. Let $\tau = \tau_n$ be defined as in (3.2), where $n$ is the total number of nulls. The goal is to show that, given $0 < \eta < 1$, as $0 < \alpha - \alpha_* \ll 1$,

$$P((1-\eta)t_* < \tau_n < (1+\eta)t_*) \to 1, \quad \text{as } n \to \infty. \tag{A.1}$$

If this holds, then it is easy to see that $G(D_{(1-\eta)t_*}) \le \text{Pow}(\alpha) \le G(D_{(1+\eta)t_*})$. As $G(D_{(1\pm\eta)t_*}) \sim (1\pm\eta)g(\mathbf{0})t_*$ and $\eta$ is arbitrary, (4.1) then follows.

The remaining part of the proof is for (A.1). By (4.2),

$$\frac{F(t)}{t} = 1 - a + \frac{aG(D_t)}{t} = \frac{1}{a_*} - Ct^{\varepsilon/K} + o(t^{\varepsilon/K}), \quad \text{as } t \to 0.$$

where $C > 0$ is a constant. By this expansion, there is $0 < \delta \ll 1$, such that

$$\inf_{s \le (1-\eta)t} \frac{F(s)}{s} > \frac{F(t)}{t}, \quad \sup_{(1+\eta)t \le s \le \delta} \frac{F(s)}{s} < \frac{F(t)}{t}, \quad \text{for } 0 < t < \delta.$$

By $g \in C(I)$ and $g(\boldsymbol{x}) < g(\mathbf{0})$ for $\boldsymbol{x} \ne \mathbf{0}$, for $t > 0$, $G(D_t)/t < g(\mathbf{0})$, yielding $F(t)/t < 1/\alpha_*$. Thus, for $0 < \alpha - \alpha_* \ll 1$, $\sup_{t \ge \delta} F(t)/t > 1/\alpha$. On the other hand, $t_* \in (0, \delta)$. Consequently,

$$\inf_{s \le (1-\eta)t_*} \frac{F(s)}{s} > \frac{1}{\alpha}, \quad \sup_{(1+\eta)t_* \le s \le 1} \frac{F(s)}{s} < \frac{1}{\alpha}.$$

Because the empirical distribution of $h(Z_i)$ converges to $F$ in probability, the above inequalities imply (A.1). □

**Corollary A.2.1.** *For procedure (3.4) using $\Gamma_u(\boldsymbol{\nu})$, $H_i$ is rejected if and only if $\boldsymbol{\nu}' \boldsymbol{\xi}_i^\varepsilon \le \zeta$, where $\zeta = \zeta_n$ is a random variable such that given $0 < \eta \ll 1$, for $0 < \alpha - \alpha_* \ll 1$, $P(|\zeta - v_*| \le \eta v_*) \to 1$ as $n \to \infty$, where*

$$v_* \sim C(\boldsymbol{\nu})(\alpha - \alpha*), \quad \text{with } C(\boldsymbol{\nu}) = \frac{K+\varepsilon}{a\alpha_*^2 g(\mathbf{0})} \bigg/ \sum \frac{\gamma_k}{\nu_k}. \tag{A.2}$$

*Proof.* Let $\zeta = h^*(\tau)$, where $h(u) = h(u; \boldsymbol{\nu})$ and $\tau$ is as in the proof of Proposition 4.1. Then $H_i$ is rejected $\iff h(\boldsymbol{\nu}'\boldsymbol{\xi}_i^\varepsilon) \leq \tau \iff \boldsymbol{\nu}'\boldsymbol{\xi}_i^\varepsilon \leq \zeta$. Let $v_* = h^*(t_*)$. Since $h^*(t) = \bar{\nu}(t/V_\varepsilon)^{\varepsilon/K}$ for $0 < t \ll 1$, the result follows from the asymptotics of $\tau$ and $t_*$. □

*Proof of Proposition 4.2.* For ease of notation, integration over $\boldsymbol{x}$ will be implicitly restricted in $I$. First consider $\text{Pow}_o(\alpha)$ of procedure (3.2) using $D_t = \Gamma_{h^*(t)}$, where $\Gamma_u = \{\boldsymbol{x} \in I : g(\boldsymbol{x}) \geq u\}$ and $h(u) = \ell(\Gamma_u)$. By Lemma A.1.1, $G(D_t)$ is strictly concave. Then by [12], $\text{Pow}_o(\alpha) = G(\Gamma_{u_*})$, where $u_* = h(t_*)$, with $t_*$ the unique positive solution to $(1-a)t + aG(D_t) = t/\alpha$. Therefore,

$$(1-a)\ell(\Gamma_{u_*}) + aG(\Gamma_{u_*}) = \ell(\Gamma_{u_*})/\alpha. \tag{A.3}$$

Using (A.3) followed by $\alpha_* = 1/(1 - a + ag(\boldsymbol{0}))$,

$$ag(\boldsymbol{0}) \int_{\Gamma_{u_*}} [\boldsymbol{\gamma}'\boldsymbol{x}^\varepsilon - r(\boldsymbol{x})] \, d\boldsymbol{x} = [1 - a + ag(\boldsymbol{0}) - 1/\alpha] \int_{\Gamma_{u_*}} d\boldsymbol{x} \tag{A.4}$$

$$= \frac{\alpha - \alpha_*}{\alpha \alpha_*} \int_{\Gamma_{u_*}} d\boldsymbol{x}. \tag{A.5}$$

Fix $\eta > 0$. For $\alpha > \alpha_*$, denote

$$v = v(\alpha) = 1 - u_*/g(\boldsymbol{0}), \quad \lambda = \alpha - \alpha_*. \tag{A.6}$$

Then $\Gamma_{u_*} = \{\boldsymbol{x} \in I : \boldsymbol{\gamma}'\boldsymbol{x}^\varepsilon - r(\boldsymbol{x}) \leq v\}$ and $v \downarrow 0$ as $\alpha \downarrow \alpha_*$. By $\gamma_k > 0$ and $r(\boldsymbol{x}) = o(|\boldsymbol{x}|^\varepsilon)$, for $0 < \alpha - \alpha_* \ll 1$, $\boldsymbol{\gamma}'\boldsymbol{x}^\varepsilon \leq (1+\eta)v$ and $|r(\boldsymbol{x})| \leq \eta\boldsymbol{\gamma}'\boldsymbol{x}^\varepsilon$ on $\Gamma_{u_*}$. Together with (A.4) and (A.5), this gives

$$(1+\eta)ag(\boldsymbol{0}) \int_{\boldsymbol{\gamma}'\boldsymbol{x}^\varepsilon \leq (1+\eta)v} \boldsymbol{\gamma}'\boldsymbol{x}^\varepsilon \, d\boldsymbol{x} \geq \frac{\lambda}{\alpha \alpha_*} \int_{\boldsymbol{\gamma}'\boldsymbol{x}^\varepsilon \leq (1-\eta)v} d\boldsymbol{x},$$

$$(1-\eta)ag(\boldsymbol{0}) \int_{\boldsymbol{\gamma}'\boldsymbol{x}^\varepsilon \leq (1-\eta)v} \boldsymbol{\gamma}'\boldsymbol{x}^\varepsilon \, d\boldsymbol{x} \leq \frac{\lambda}{\alpha \alpha_*} \int_{\boldsymbol{\gamma}'\boldsymbol{x}^\varepsilon \leq (1+\eta)v} d\boldsymbol{x}.$$

By Lemma 3.2, the inequalities imply

$$\frac{a(1+\eta)^{K/\varepsilon+2}Kv}{K+\varepsilon} \geq \frac{\lambda(1-\eta)^{K/\varepsilon}}{\alpha_* \alpha g(\boldsymbol{0})}, \quad \frac{a(1-\eta)^{K/\varepsilon+2}Kv}{K+\varepsilon} \leq \frac{\lambda(1+\eta)^{K/\varepsilon}}{\alpha_* \alpha g(\boldsymbol{0})}.$$

Since $\eta$ is arbitrary, it follows that

$$v \sim \frac{(K+\varepsilon)\lambda}{Ka\alpha_*^2 g(\boldsymbol{0})}, \quad \text{as } \alpha \downarrow \alpha_*.$$

Comparing with (A.2), $v \sim C(\boldsymbol{\gamma})\lambda$.

Applying (A.3) followed by (3.10),

$$\text{Pow}_o(\alpha) = \left(\frac{1}{a\alpha} - \frac{1}{a} + 1\right) \int_{\boldsymbol{\gamma}'\boldsymbol{x}^\varepsilon - r(\boldsymbol{x}) \leq v} d\boldsymbol{x}$$

$$\sim g(\boldsymbol{0}) \int_{\boldsymbol{\gamma}'\boldsymbol{x}^\varepsilon \leq C(\boldsymbol{\gamma})\lambda} d\boldsymbol{x} \sim g(\boldsymbol{0})V_\varepsilon \left[\frac{1}{\bar{\nu}} \frac{(K+\varepsilon)\lambda}{Ka\alpha_*^2 g(\boldsymbol{0})}\right]^{K/\varepsilon}.$$

On the other hand, by Proposition 4.1, $\text{Pow}(\alpha)$ has the same asymptotic. Then $\text{Pow}_o(\alpha)/\text{Pow}(\alpha) \to 1$.

Let $\vartheta = C(\boldsymbol{\gamma})(\alpha - \alpha_*)$. The above argument shows that, for procedure (3.2) based on $\Gamma_u$, if $0 < \alpha - \alpha_* \ll 1$, then $H_i$ is rejected if only if $\boldsymbol{\gamma}'\boldsymbol{\xi}_i^\varepsilon - r(\boldsymbol{\xi}_i) \leq \zeta_o$, where $\zeta_o = \zeta_{o,n}$ is a random variable satisfying $P(|\zeta_o - \vartheta| \leq \eta\vartheta) \to 1$ as $n \to \infty$. On the other hand, by Corollary A.2.1, for procedure (3.9) based on $\Gamma_u(\boldsymbol{\gamma})$, $H_i$ is rejected if and only if $\boldsymbol{\gamma}'\boldsymbol{\xi}_i^\varepsilon \leq \zeta$, where $\zeta$ satisfies $P(|\zeta - \vartheta| \leq \eta\vartheta\}) \to 1$ as $n \to \infty$.

Recall $\mathcal{V}_o = \{\text{true } H_i : \boldsymbol{\gamma}'\boldsymbol{\xi}_i^\varepsilon - r(\boldsymbol{\xi}_i) \leq \zeta_o\}$ and $\mathcal{V} = \{\text{true } H_i : \boldsymbol{\gamma}'\boldsymbol{\xi}_i^\varepsilon \leq \zeta\}$. Since $r(\boldsymbol{x}) = o(\boldsymbol{\gamma}'\boldsymbol{x}^\varepsilon)$, by the asymptotics of $\zeta_o$ and $\zeta$, given $0 < \eta \ll 1$, as $0 < \alpha - \alpha_* \ll 1$,

$$P(\{\text{true } H_i : \boldsymbol{\gamma}'\boldsymbol{\xi}_i^\varepsilon \leq (1-\eta)\vartheta\} \subset \mathcal{V}_o \cap \mathcal{V}) \to 1$$
$$P(\mathcal{V}_o \triangle \mathcal{V} \subset \{\text{true } H_i : (1-\eta)\vartheta < \boldsymbol{\gamma}'\boldsymbol{\xi}_i^\varepsilon \leq (1+\eta)\vartheta\}) \to 1.$$

It follows that

$$r_V(\alpha) \leq \int_{(1-\eta)\vartheta \leq \boldsymbol{\gamma}'\boldsymbol{x}^\varepsilon \leq (1+\eta)\vartheta} d\boldsymbol{x} \bigg/ \int_{\boldsymbol{\gamma}'\boldsymbol{x}^\varepsilon \leq (1-\eta)\vartheta} d\boldsymbol{x} = 1 - \left(\frac{1+\eta}{1-\eta}\right)^{K/\varepsilon}.$$

As $\eta$ is arbitrary, this gives $r_V(\alpha) \to 0$ as $\alpha \downarrow \alpha_*$. Likewise, $r_D(\alpha) \to 0$. □

*Proof of Lemma 4.1.* By the assumptions of the lemma and [12],

$$\text{Pow}_i(\alpha) = G(D_{it_i^*}), \tag{A.7}$$
$$\text{where} \quad t_i^* = \sup\{t : (1-a)t + aG(D_{it}) \geq t/\alpha\}$$

and furthermore, $t_i^* < T$. Since $G(D_{1t}) < G(D_{2t})$ for $t < T$ and both are continuous, it is seen that $t_1^* < t_2^*$. On the other hand,

$$G(D_{it_i^*}) = \left(\frac{1}{a\alpha} - \frac{1}{a} + 1\right) t_i^* \sim g(\boldsymbol{0}) t_i^*, \quad \text{as} \quad \alpha \to \alpha_*. \tag{A.8}$$

As a result, $\text{Pow}_1(\alpha) < \text{Pow}_2(\alpha)$. □

*Proof of Proposition 4.3.* By (A.8),

$$\frac{t_1^*}{t_2^*} = \frac{G(D_{1t_1^*})}{G(D_{2t_2^*})} = \frac{\text{Pow}_1(\alpha)}{\text{Pow}_2(\alpha)} \quad \Longrightarrow \quad \frac{t_1^*}{t_2^*} = \frac{g(0)t_1^* - G(D_{1t_1^*})}{g(0)t_2^* - G(D_{2t_2^*})}.$$

Suppose $M < \infty$. Then by the last equality and (4.2), as $\alpha \downarrow \alpha_*$,

$$\frac{t_1^*}{t_2^*} \sim M \times \frac{g(0)t_1^* - G(D_{2t_1^*})}{g(0)t_2^* - G(D_{2t_2^*})} \sim M \left(\frac{t_1^*}{t_2^*}\right)^{\varepsilon/K+1},$$

giving $t_1^*/t_2^* \sim (1/M)^{K/\varepsilon}$ and the proof. The case $M = \infty$ is likewise. □

*Proof for Example 4.1* (Case $K > 2$). Define $D_{1t}$ and $D_{2t}$ as in the case $K = 2$. Fix $\eta > 0$ and $a_k \in (0, \gamma_k)$, so that $g(\boldsymbol{x}) < g(\boldsymbol{0})(1 - \sum a_k x_k^\varepsilon)$ on $[0, \eta]^K$. Then

$$G(D_{1t}) \leq g(\boldsymbol{0}) \left[ t - \sum a_k \int_{\substack{x_1 \cdots x_K \leq s \\ \text{all } x_i \in (0, \eta)}} x_k^\varepsilon \, d\boldsymbol{x} \right]$$

where $s = \exp(-F_K^{-1}(1-t))$, with $F_K$ the Gamma distribution function with $K$ degrees of freedom and scale parameter 1. We need to evaluate

$$\int_{\substack{x_1 \cdots x_K \leq s \\ \text{all } x_i \in (0, \eta)}} x_k^\varepsilon \, d\boldsymbol{x} = \eta^{K+\varepsilon} \int_{\substack{x_1 \cdots x_K \leq r \\ \text{all } x_i \in (0, 1)}} x_k^\varepsilon \, d\boldsymbol{x} = \eta^{K+\varepsilon} E[U^\varepsilon \mathbf{1}\{UV \leq r\}],$$

where $r = s/\eta^K$, $U \sim \text{Unif}(0, 1)$ and $V$ are independent, and $V$ is the product of $U_1, \ldots, U_{K-1}$ i.i.d. $\sim \text{Unif}(0, 1)$. By transformation $X = -\ln U$ and $Z = -\ln V$, the expectation equals

$$A(r) := \int_0^{-\ln r} e^{-(1+\varepsilon)(-\ln r - x)}[1 - F_{K-1}(x)] \, dx + \int_{-\ln r}^\infty e^{-(1+\varepsilon)x} \, dx$$

Recall that for $n \geq 1$, as $x \to \infty$,

$$1 - F_n(x) \sim x^{n-1} e^{-x}/(n-1)! \tag{A.9}$$

As $t \to 0$, $s \to 0$ and hence $r \to 0$, yielding

$$A(r) \sim \frac{r^{1+\varepsilon}}{(K-2)!} \int_0^{-\ln r} x^{K-2} e^{\varepsilon x} \, dx + r^{1+\varepsilon} \sim \frac{r(\ln r^{-1})^{K-2}}{(K-2)!}.$$

So for $t \ll 1$, $G(D_{1t}) \leq g(\boldsymbol{0})[t - Cs(\ln s^{-1})^{K-2}]$, where $C > 0$ is a constant. On the other hand, by (4.2), $G(D_{2t}) \geq g(\boldsymbol{0})[t - C' t^{1+\varepsilon/K}]$, where $C' > 0$ is another constant. By (A.9), for some constants $c_1, c_2 > 0$,

$$\frac{s(\ln s^{-1})^{K-2}}{t^{1+\varepsilon/K}} = \frac{s(\ln s^{-1})^{K-2}}{(1 - F_K(-\ln s))^{1+\varepsilon/K}} \sim \frac{c_1}{s^{\varepsilon/K}(-\ln s)^{c_2}} \to \infty.$$

As a result, for $t \ll 1$, $G(D_{1t}) < G(D_{2t})$. $\square$

*Proof for Example 4.2.* Recall $\bar{\Phi}^{-1}(s) \sim \sqrt{2 \log(1/s)}$ as $s \downarrow 0$. Let $s_k = w_k/w$. Given $c > 1$ and $a_k \in (0, \gamma_k)$, there is $0 < \eta \ll 1$, such that 1) $g(\boldsymbol{x}) < g(\boldsymbol{0})(1 - \boldsymbol{a}'\boldsymbol{x})$ for $\boldsymbol{x} \in (0, \eta)^K$ and 2) letting

$$B_t := \left\{ \boldsymbol{x} \in (0, \eta)^K : \sum_k s_k \sqrt{\log(1/x_k)} \geq c\sqrt{\log(1/t)} \right\},$$

for $0 < t \ll 1$, $B_t \subset (0, \eta)^K \cap D_{1t}$. Then

$$G(D_{1t}) \leq g(\boldsymbol{0}) \left[ \int_{D_{1t}} d\boldsymbol{x} - \int_{B_t} \boldsymbol{a}'\boldsymbol{x} \, d\boldsymbol{x} \right] = g(\boldsymbol{0}) \left[ t - \sum_{k=1}^K a_k \int_{B_t} x_k^\varepsilon \, d\boldsymbol{x} \right].$$

Let $x_k = t^{y_k}$ and $\delta_t = \ln \eta / \ln t$. For each $j = 1, \ldots, K$,

$$\int_{B_t} x_j^\varepsilon \, d\boldsymbol{x} = |\ln t|^K \int_{\substack{y_k \in (\delta_t, \infty) \\ \sum_k s_k \sqrt{y_k} \geq c}} e^{(y_1 + \cdots + y_K + \varepsilon y_j) \ln t} \, d\boldsymbol{y}.$$

As $t \downarrow 0$, $\ln t \to -\infty$ and $\delta_t \downarrow 0$. So by $\sum s_i^2 = 1$,

$$\frac{1}{\ln t} \ln \int_{B_t} x_j^\varepsilon \, d\boldsymbol{x} \to \inf_{\sum_k s_k \sqrt{y_k} \geq c} [y_1 + \cdots + y_K + \varepsilon y_j] = c^2 \left[1 - \frac{\varepsilon s_j^2}{1 + \varepsilon}\right]^{-1}$$

and hence

$$\varlimsup_{t \downarrow 0} \frac{\ln[g(\mathbf{0})t - G(D_{1t})]}{\ln t} \leq \lim_{t \downarrow 0} \frac{1}{\ln t} \ln \left(\sum_k a_k \int_{B_t} x_k^\varepsilon \, d\boldsymbol{x}\right)$$

$$= \min_k \lim_{t \downarrow 0} \frac{1}{\ln t} \ln \left(\int_{B_t} x_k^\varepsilon \, d\boldsymbol{x}\right) = c^2 \left[1 - \frac{\varepsilon}{1 + \varepsilon} \min_k s_k^2\right]^{-1}.$$

Since $c > 1$ is arbitrary and $\min_k s_k^2 \leq 1/K$, (4.3) then follows. □

*Proof of Equations (4.5) and (4.7).* By (4.4), as $t \downarrow 0$,

$$G(D_{1t}) = g(\mathbf{0}) \int_{0 \leq x_k \leq f_k(t)} (1 - \boldsymbol{\gamma}' \boldsymbol{x}^\varepsilon + r(\boldsymbol{x})) \, d\boldsymbol{x}$$

$$= g(\mathbf{0}) \left(t - \frac{1 + o(1)}{1 + \varepsilon} \sum \gamma_k f_k(t)^{1+\varepsilon} \prod_{i \neq k} f_i(t)\right).$$

Since $\prod f_k(t) = t$, (4.5) then follows.

Let $f_k(t) = c_k t^{1/K}$. By (A.7), (A.8) and (4.5), $\text{Pow}(\alpha) \sim g(\mathbf{0}) t^*$ as $\alpha \to \alpha_*$, where $t^* > 0$ is the solution to

$$1 - a + ag(\mathbf{0}) \left(1 - \frac{1}{1 + \varepsilon} \sum \gamma_k c_k^\varepsilon t^{\varepsilon/K}\right) = 1/\alpha.$$

Since $1/\alpha_* = 1 - a + ag(\mathbf{0})$,

$$t_*^{\varepsilon/K} = \frac{1 + \varepsilon}{ag(\mathbf{0})} \left(\frac{1}{\alpha_*} - \frac{1}{\alpha}\right) \bigg/ \sum \gamma_k c_k^\varepsilon \sim \left(\frac{1 + \varepsilon}{a \alpha_*^2 g(\mathbf{0})} \bigg/ \sum \gamma_k c_k^\varepsilon\right)(\alpha - \alpha_*).$$

It is then easy to see (4.7) holds. □

*Proof for Example 4.4.* We need to show $L < g(\mathbf{0}) = \prod_k \rho_k$, where

$$L = \sup_x \frac{\sum_k f_k(x) \prod_{j \neq k} F_j(x)}{\sum_k f_{0k}(x) \prod_{j \neq k} F_{0j}(x)} \leq \max_k \sup_x \frac{f_k(x) \prod_{j \neq k} F_j(x)}{f_{0k}(x) \prod_{j \neq k} F_{0j}(x)}$$

Let $\psi_k = F_k/F_{0k}$. By $F_k(x) = \int_{-\infty}^x r_k(t)f_{0k}(t)\,dt$ and the assumption (4.8) on $r_k$, for any finite $c$, $\sup_{x \le c} \psi_k(x) < \rho_k$. Since $\psi_k(x)$ is continuous and tends to 1 as $x \to \infty$, it follows that $\sup_x \psi_k(x) < \rho_k$. Then for each $k$,

$$\sup_x \frac{f_k(x)\prod_{j\ne k} F_j(x)}{f_{0k}(x)\prod_{j\ne k} F_{0j}(x)} = \sup_x r_k(x) \prod_{j\ne k} \psi_j(x) < \prod \rho_k = g(\mathbf{0}).$$

Thus $L < g(\mathbf{0})$. □

*Proof for Example 4.5.* We need to show $L := \sup_{x \ge T} r(x) < g(\mathbf{0}) = \rho_1\rho_2$ for $T > 0$, where

$$r(x) := \frac{\int_0^x f_1(t)f_2(x-t)\,dt}{\int_0^x f_{01}(t)f_{02}(x-t)\,dt},$$

As $f_k(x) < \rho_k f_{0k}(x)$ for any $x > 0$ with $f_k(x) > 0$, it is seen $r(x) < g(\mathbf{0})$ for any $x > 0$. It remains to be shown $r(\infty) < g(\mathbf{0})$. Once this is done, by the continuity of $r(x)$ on $(0, \infty)$, $L < g(\mathbf{0})$. Given $c, \lambda \in (0, 1)$,

$$I := \int_0^x f_{01}(t)f_{02}(x-t)\,dt = I_1 + I_2$$

where $I_1$ is the integral over $[0, cx]$ and $I_2$ over $[cx, x]$. As $x \gg 1$, for $t \in [0, cx]$, $\lambda x^{-s_2} \le f_{02}(x-t) \le \lambda^{-1}(1-c)^{-s_2} x^{-s_2}$. Then

$$\lambda x^{-s_2} \int_0^{cx} f_{01}(t)\,dt \le I_1 \le \lambda^{-1}(1-c)^{-s_2} x^{-s_2} \int_0^{cx} f_{01}(t)\,dt.$$

Similarly, by using $f_{01}(x) \sim x^{-s_1}$,

$$\lambda x^{-s_1} \int_{cx}^x f_{02}(x-t)\,dt \le I_2 \le \lambda^{-1} c^{-s_1} x^{-s_1} \int_{cx}^x f_{02}(x-t)\,dt$$

As $x \to \infty$, $\int_0^{cx} f_{01} \to 1$ and $\int_{cx}^x f_{02}(x-t)\,dt = \int_0^{(1-c)x} f_{02} \to 1$. Therefore, if $s_1 > s_2$, then $I_2 = o(I_1)$ and hence $I \sim I_1$. Since $c$ and $\lambda$ are arbitrary, $I \sim x^{-s_2}$. Likewise

$$J := \int_0^x f_1(t)f_2(x-t)\,dt \sim \rho_2 x^{-s_2},$$

and hence $r(x) \to \rho_2$. Similarly, if $s_1 < s_2$, then $r(x) \to \rho_1$. If $s_1 = s_2$, write $I = I_1 + I_2 + I_3$, where the integrals $I_i$ are over $[0, cx]$, $[cx, (1-c)x]$ and $[(1-c)x, x]$ respectively, with $0 < c \ll 1$. Following argument similar to the above, $I_1 \sim I_3 \sim x^{-s_1}$ while $I_2 = O(x^{-2s_1})$. Then $I \sim 2x^{-s_1}$. Likewise, $J \sim (\rho_1 + \rho_2)x^{-s_1}$. Thus $r(x) \to (\rho_1 + \rho_2)/2$. In any case, $r(\infty) < g(\mathbf{0})$. □

### A.3. Some facts about $V_\varepsilon$

Let $V_{\varepsilon,K}(u) = \int_0^1 \cdots \int_0^1 \mathbf{1}\left\{\sum x_k^\varepsilon \leq u\right\} d\boldsymbol{x}$ and $V_{\varepsilon,K} = V_{\varepsilon,K}(1)$. It is not hard to see that $V_{\varepsilon,K}(u) = u^{K/\varepsilon} V_{\varepsilon,K}$ for $u \in [0,1]$. It can be shown that

$$V_{\varepsilon,K} = \frac{(1/\varepsilon)^{K-1}\Gamma(1/\varepsilon)^K}{K\Gamma(K/\varepsilon)}, \quad K = 1, 2, \ldots. \tag{A.10}$$

This is clear for $K = 1$. For $K > 1$, by first integrating $x_1, \ldots, x_{K-1}$,

$$V_{\varepsilon,K} = \int_0^1 V_{\varepsilon,K-1}(1 - x_K^\varepsilon) dx_K = V_{\varepsilon,K-1} \int_0^1 (1 - x^\varepsilon)^{(K-1)/\varepsilon} dx$$

$$= V_{\varepsilon,K-1} \times \frac{1}{\varepsilon} \int_0^1 t^{1/\varepsilon - 1}(1-t)^{(K-1)/\varepsilon} dt$$

$$= V_{\varepsilon,K-1} \times \frac{(K-1)\Gamma(1/\varepsilon)\Gamma((K-1)/\varepsilon)}{\varepsilon K\Gamma(K/\varepsilon)}.$$

Then (A.10) follows by induction.

Finally, in Example 4.3, it is claimed that $\frac{(1/V_{\varepsilon,K})^{\varepsilon/K}}{K+\varepsilon} < \frac{1}{1+\varepsilon}$ for $K > 1$. To show this, write $t = 1/\varepsilon$ and $H(t) = \ln[t\Gamma(t)/(1+t)^t]$, Then by (A.10), the above inequality is equivalent to $KH(t) \geq H(Kt)$, $K > 1$. It is not hard to get $\lim_{t\to 0+} H(t) = 0$. Therefore, if one can show $H(t)$ is concave for $t > 0$, then the desired inequality is obtained. Now

$$H''(t) = (\ln \Gamma(t))'' - \frac{1}{t^2} - \frac{1}{(1+t)^2} - \frac{1}{1+t}.$$

It is known that $(\ln \Gamma(t))'' = \sum_{k=0}^\infty (k+t)^{-2}$, $t > 0$. Since $1/(k+t)^2 < 1/(k-1+t) - 1/(k+t)$ for $k \geq 2$, then it is seen that $H''(t) < 0$ for $t > 0$ and hence $H(t)$ is strictly concave.